\documentclass[12pt, reqno]{amsart}
\usepackage{amssymb,amsthm,amsmath}
\usepackage[numbers,sort&compress]{natbib}
\usepackage{amssymb,amsmath}
\usepackage{amsfonts}
\usepackage{mathrsfs}
\usepackage{latexsym}
\usepackage{amssymb}
\usepackage{amsthm}
\usepackage{indentfirst}
\date{April 5, 2016}

\let\oldsection\section
\renewcommand\section{\setcounter{equation}{0}\oldsection}

\newtheorem{corollary}{Corollary}[section]
\newtheorem{theorem}{Theorem}[section]
\newtheorem{lemma}{Lemma}[section]
\newtheorem{proposition}{Proposition}[section]
\newtheorem{definition}{Definition}[section]

\newtheorem{example}{Example}[section]
\newtheorem{remark}{Remark}[section]

\allowdisplaybreaks
\begin{document}

\title[Recent Advances Concerning Certain Class of Geophysical Flows]{Recent advances concerning certain class of geophysical flows}


\author{Jinkai~Li}
\address[Jinkai~Li]{Department of Computer Science and Applied Mathematics, Weizmann Institute of Science, Rehovot 76100, Israel.}
\email{jklimath@gmail.com}

\author{Edriss~S.~Titi}
\address[Edriss~S.~Titi]{
Department of Mathematics, Texas A\&M University, 3368 TAMU, College Station, TX 77843-3368, USA. ALSO, Department of Computer Science and Applied Mathematics, Weizmann Institute of Science, Rehovot 76100, Israel.}
\email{titi@math.tamu.edu and edriss.titi@weizmann.ac.il}

\keywords{Anisotropic incompressible Navier-Stokes equations; primitive equations for oceanic and atmospheric flow; geophysical flow; Boussinesq equations; tropical atmosphere; hydrostatic approximation; singular perturbation limit.}
\subjclass[2010]{35A01,
35B45, 35Q30, 35Q86, 76D03, 76D09.}


\begin{abstract}
This paper is devoted to reviewing several recent developments
concerning certain class of geophysical models,
including the
primitive equations (PEs) of atmospheric and oceanic dynamics and a tropical atmosphere model.
The PEs for large-scale oceanic and atmospheric dynamics
are derived from the
Navier-Stokes equations coupled to the heat convection by adopting the Boussinesq and hydrostatic
approximations, while the tropical atmosphere model considered here
is a nonlinear interaction system between the barotropic mode and the
first baroclinic mode of the tropical atmosphere with moisture.

We are mainly concerned with the global well-posedness of strong solutions
to these systems, with full or partial viscosity, as well as certain singular perturbation small parameter limits related to these systems,
including the small aspect ratio limit from
the Navier-Stokes equations to the PEs, and a small relaxation-parameter in the tropical atmosphere model. These limits provide a rigorous justification to the hydrostatic balance in the PEs, and to the relaxation limit of
the tropical atmosphere model, respectively. Some conditional uniqueness of
weak solutions, and
the global well-posedness of weak solutions with certain class of
discontinuous initial data, to the PEs are also presented.
\end{abstract}

\maketitle

\tableofcontents
\allowdisplaybreaks

\section{The primitive equations (PEs) of atmospheric and oceanic dynamics}
\label{SECPEINT}
In the atmospheric and oceanic dynamics, the Boussinesq approximation
model is accepted as the fundamental model that governs their motion. This
system in the case of incompressible flows reads as
\begin{equation}
  \label{BE}
  \left\{
  \begin{array}{l}
    \partial_t\mathcal U+(\mathcal U\cdot\nabla)\mathcal U+\nabla p-\nu_h\Delta_H\mathcal U-\nu_z\partial_z^2\mathcal U+f_0k\times\mathcal U=\theta e_3,\\
    \nabla\cdot\mathcal U=0,\\
    \partial_t\theta+\mathcal U\cdot\nabla\theta-\kappa_h\Delta_H\theta -\kappa_z\partial_z^2\theta=0,
  \end{array}
  \right.
\end{equation}
in the space-time domain $\mathbb R^3\times(0,\infty)$, where the
unknowns are the velocity field $\mathcal U=(v,w)$, with horizontal
velocity $v=(v^1, v^2)$ and vertical velocity $w$, the pressure $p$,
and the temperature $\theta$. The nonnegative constants $\nu_h, \nu_z,
\kappa_h$ and $\kappa_z$ are the horizontal viscosity, vertical
viscosity, horizontal diffusivity and vertical diffusivity
coefficients, respectively, $f_0$ is
the Coriolis parameter and $k=(0,0,1)$. Note that here we have
taken the earth rotation into consideration. We denote by $(x,y,z)$
and $t$, respectively, as the spatial and temporal variables, and
use $\Delta_H=\partial_x^2+\partial_y^2$ to denote the horizontal
Laplacian.

For large-scale oceanic and atmospheric dynamics, the vertical scale ($10$\,-\,$20$ kilometers) is much smaller than the horizontal scales (several thousands of kilometers), and consequently, the aspect ratio, i.e.\,the ratio of the depth (or hight) to the horizontal width, is very small.
In view of this thinness (the ratio of the vertical hight to horizontal width is small) of the ocean and
atmosphere, we consider the above Boussinesq equations (\ref{BE}) in a
thin domain $\Omega_\varepsilon:=M\times(-\varepsilon, \varepsilon)$, with $M=(0,L_1)\times(0,L_2)$, for two positive constants $L_1$ and
$L_2$, and $\varepsilon$ a small positive parameter.
Following Az\'erad--Guill\'en
\cite{AZGU}, we suppose that the viscosity and diffusivity coefficients have orders $(\nu_h, \kappa_h)=O(1)$ and $(\nu_z, \kappa_z)=O(\varepsilon^2)$, and for
simplicity, we assume $\nu_h=\kappa_h=1$ and $\nu_z=\kappa_z=\varepsilon^2$. By rescaling, we introduce the following new unknowns
\begin{eqnarray*}
  && v_\varepsilon(x,y,z,t)=v(x,y,\varepsilon z,t), \quad w_\varepsilon(x,y,z,t)=\frac1\varepsilon w(x,y,\varepsilon z,t), \\
  &&p_\varepsilon(x,y,z,t)=p(x,y,\varepsilon z,t),\quad \theta_\varepsilon(x,y,z,t)=\varepsilon\theta(x,y,\varepsilon z,t),
\end{eqnarray*}
for $(x,y,z)\in\Omega:= M\times(-1,1)$ and $t\in(0,\infty)$, Boussinesq
system (\ref{BE}), defined on the $\varepsilon$-dependent domain $\Omega_\varepsilon$, can be transformed to the following scaled Boussinesq equations
\begin{equation}
  \label{SBE}
  \left\{
  \begin{array}{l}
    \partial_t v_\varepsilon+(v_\varepsilon \cdot\nabla_H)  v_\varepsilon+w_\varepsilon\partial_zv_\varepsilon - \Delta   v_\varepsilon +\nabla_H p_\varepsilon+f_0k\times v_\varepsilon =0,\\
    \nabla_H\cdot v_\varepsilon+\partial_zw_\varepsilon=0,\\
    \varepsilon^2(\partial_tw_\varepsilon+v_\varepsilon\cdot\nabla_Hw_\varepsilon+ w_\varepsilon\partial_zw_\varepsilon-\Delta w_\varepsilon)+\partial_zp_\varepsilon=\theta_\varepsilon,\\
    \partial_t\theta_\varepsilon+v_\varepsilon\cdot\nabla_H \theta_\varepsilon+w_\varepsilon\partial_z\theta_\varepsilon -\Delta\theta_\varepsilon=0,
  \end{array}
  \right.
\end{equation}
in the fixed space-time domain $\Omega\times(0,\infty)$. Here, the
rotating Coriolis term $f_0k\times v_\varepsilon$, in (\ref{SBE}), is understood
as the horizontal components of the corresponding term, that is $k\times v_\varepsilon=(-v_{\varepsilon}^2,
v_{\varepsilon}^1)$.

Let us take the formal limit, $\varepsilon\rightarrow0^+$, and suppose that $(v_\varepsilon, w_\varepsilon, p_\varepsilon, \theta_\varepsilon)$ converges to $(v, w, p, \theta)$ in a suitable
sense, then the vertical momentum equation in (\ref{SBE}) degenerates to the following hydrostatic balance
\begin{equation}
  \label{HB}
\partial_zp=\theta,
\end{equation}
and as a result, one obtains the following system, known as the primitive equations (PEs)
\begin{equation}
  \label{PES}
  \left\{
  \begin{array}{l}
    \partial_tv +(v \cdot\nabla_H)v+w\partial_zv-\Delta v+\nabla_H
    p+f_0k\times v=0,\\
    \nabla_H\cdot v+\partial_zw=0,\\
    \partial_zp=\theta,\\
    \partial_t\theta+v\cdot\nabla_H\theta+w\partial_z\theta
    -\Delta\theta=0,
  \end{array}
  \right.
\end{equation}
in $\Omega\times(0,\infty)$, where as before, the rotating term is understood as $k\times v=(-v_2, v_1)$.

The above small aspect ratio limit,
from (\ref{SBE}) to (\ref{PES}), can be rigorously justified.
Indeed, the weak convergence of such limit was shown in Az\'erad--Guill\'en, \cite{AZGU}, while the strong convergence was proved by Li--Titi \cite{LITITIHYDRO}.
This strong convergence result, reported in \cite{LITITIHYDRO}, will be presented in section
\ref{SECPEHYD}. We remark that it is
necessary to consider the anisotropic viscosities, so that the
Navier-Stokes equations converge to the PEs,
as the aspect ratio $\varepsilon$ goes to zero.
In fact, for the isotropic case that $(\mu,\nu)=O(1)$, it has been
shown in Bresh--Lemoine--Simon \cite{BRLESI} that the stationary
Navier-Stokes equations
converge to a linear system with only vertical dissipation, instead of
the stationary PEs.

As a result of the above discussion, the primitive equations form a fundamental block in models of the oceanic and atmospheric
dynamics, see, e.g., the books Haltiner--Williams \cite{HAWI},
Lewandowski \cite{LEWAN}, Majda \cite{MAJBOOK}, Pedlosky \cite{PED},
Vallis \cite{VALLIS}, Washington--Parkinson \cite{WP} and
Zeng \cite{ZENG}. The mathematical studies of the PEs were started
by Lions--Temam--Wang \cite{LTW92A,LTW92B,LTW95}
in the 1990s, where among other issues, global existence of
weak solutions was established; however, the uniqueness of
weak solutions is still an open question, even
for the two-dimensional case. This is different from the
incompressible Navier-Stokes equations, as it is well-known that
the weak solutions to the two-dimensional incompressible
Navier-Stokes equations are unique (see, e.g., Constantin--Foias \cite{CONFOINSBOOK}, Ladyzhenskaya \cite{LADYZHENSKAYA} and Temam \cite{TEMNSBOOK}). Moreover, it has been shown by Bardos et al.\,\cite{BLNNT} that if the initial
data of the three-dimensional Navier-Stokes equations is a function
of only two spatial variables, then the Leray-Hopf weak solution is unique and
remains a function of only two spatial variables.
The main obstacle of proving
the uniqueness of weak solutions to the PEs is the absence of the dynamical equation for
the vertical velocity. In fact, the vertical velocity can only
be recovered from the horizontal velocity through the
incompressibility condition, and as a result, there is one
derivative loss for the horizontal velocity.
Though the general result on the uniqueness of weak solutions to
the PEs is still unknown, some particular cases
have been solved, see Bresch et al.\,\cite{BGMR03}, Petcu--Temam--Ziane \cite{PTZ09} and Tachim Madjo \cite{TACHIM} for the case of the
so-called $z$-weak solutions, i.e., the weak solutions with initial
data in $X=\{v_0\in L^6|\partial_zv_0\in L^2\}$,
and Kukavica--Pei--Rusin--Ziane \cite{KPRZ} for the case of
weak solutions with continuous initial data. In a recent work \cite{LITITIUNIDIS}, we generalize the above results and show that weak solutions to the PEs, with initial data taken as small $L^\infty$ perturbations of functions in $X$, are unique, and this result will be presented in section \ref{SECPEUNIDIS}.

In the context of strong solutions, the local well-posedness
was established in Guill\'en-Gonz\'alez et al \cite{GMR01} (see also Temam--Ziane \cite{RZ04}). Remarkably, taking advantage of the fact that effectively the unknown pressure is a function of only two spatial variables, it was first shown in Cao--Titi \cite{CAOTITI07} the global regularity of strong solutions, with full viscosity and diffusivity, of the 3D PEs with the relevant physical boundary conditions (see, also,
Kobelkov \cite{KOB06}). Taking advantage of the above observation, made in \cite{CAOTITI07} concerning the pressure, it was shown by Kukavica--Ziane \cite{KZ07A,KZ07B} the global regularity of the 3D PEs subject to Dirichlet boundary conditions. Moreover, Hieber--Kashiwabara \cite{HIEKAS} have made some recent progress towards relaxing the smoothness on the initial
data, but still for the system with full viscosity and diffusivity, of the PEs by using the semigroup method. As it has been mentioned above, the key observation made in \cite{CAOTITI07} for the
global existence of strong solutions to the PEs in 3D is that effectively
the pressure,
thanks to the hydrostatic balance (\ref{HB}), is a function of only two spatial variables, i.e.\,depends only on the horizontal spatial variables and the time variable,
up to
a term determined by the temperature; note that the temperature
satisfies the maximal principle, therefore the pressure for the PEs has
better properties than that for the Navier-Stokes equations.

In the oceanic and atmospheric dynamics, due to the strong horizontal turbulent mixing,
the horizontal viscosity and diffusivity are much stronger than the vertical viscosity and diffusivity, respectively.
Therefore, both physically and mathematically, it is important to study the limiting case that the vertical viscosity or diffusivity vanishes.
In view of this, some further developments on the global existence of strong solutions
to the
PEs, with full viscosity but partial diffusivity, have been made by
Cao--Titi \cite{CAOTITI12} and Cao--Li--Titi
\cite{CAOLITITI1,CAOLITITI2}.
Furthermore, the recent works by Cao--Li--Titi
\cite{CAOLITITI3,CAOLITITI4,CAOLITITI5} show
that the horizontal viscosity turns out
to be more crucial than the vertical one for guaranteeing
the global well-posedness,
as a results they show that
merely the horizontal viscosity is sufficient to guarantee the global
well-posedness of strong solutions to the PEs, as long as one still has
either the horizontal diffusivity or the vertical diffusivity.
These results on the global well-posedness of strong solutions to
the PEs, with full or only horizontal viscosities and partial diffusivity, will be presented in section \ref{SECPEFP} and section \ref{SECPEHP}.

Notably, smooth solutions to
the inviscid PEs, with or without coupling to the temperature equation, have
been shown by Cao et al.\,\cite{CINT} and Wong \cite{WONG} to blow up in
finite time. For the results on the local well-posedness with monotone
or analytic initial data of the
inviscid PEs (which are also called inviscid Prandtl equations or
hydrostatic Euler equations), i.e.,
system (\ref{PES}), without the Lapalacian term $\Delta v$, see,
e.g., Brenier \cite{BRE}, Masmoudi--Wong \cite{MASWON},
Kukavica--Temam--Vicol--Ziane \cite{KTVZ}, and the references
therein.

Another
remarkable difference between the incompressible Navier-Stokes equations
and the PEs is their well-posedness theories with $L^p$
initial data. It is well-known that, for any $L^p$ initial data, with
$d\leq p\leq\infty$, there is a unique local mild solution to the
$d$-dimensional incompressible Navier-Stokes equations, see Kato \cite{KATO} and Giga \cite{GIGA86,GIGA99};
however, for the PEs, though the $L^p$ norms of the
weak solutions remain finite, up to any finite time, as long as the
initial data belong to $L^p$ (one can apply Proposition 3.1 in \cite{CAOLITITI3}, see Proposition \ref{PROPPEHH} in section \ref{SECPEHH}, below, to achieve this fact, with the help of some regularization procedure), it is
still an open question to show the uniqueness of weak solutions to
the PEs with $L^p$ initial data. However, for the case that $p=\infty$, as we mentioned above,
it has been shown in \cite{LITITIUNIDIS} that
the weak solutions to the PEs, with initial data in $L^\infty$, are unique, as long as the discontinuity of the initial data is sufficiently small.

\section{The PEs with full viscosity but partial diffusivity}
\label{SECPEFP}
This section and the next one are devoted to the study of global well-posedness of strong solutions to the PEs, with full or only horizontal viscosity and partial diffusivity. In this section, we consider the full viscosity case (however, with only partial diffusivity), and the case with only horizontal viscosity is postponed
to the next section, section \ref{SECPEHP}.

\subsection{The vertical diffusivity case with full viscosity}
\label{SECPEFV}
Consider the following version of the PEs with full viscosity, but only vertical diffusivity
\begin{equation}
  \label{PEFV'}
  \left\{
  \begin{array}{l}
    \partial_tv +(v \cdot\nabla_H)v+w\partial_zv-\Delta v+\nabla_H
    p+f_0k\times v=0,\\
    \nabla_H\cdot v+\partial_zw=0,\\
    \partial_zp=T,\\
    \partial_tT+v\cdot\nabla_HT+w\partial_zT
    -\partial_z^2T=0,
  \end{array}
  \right.
\end{equation}
in $\Omega_0=M\times(-h,0)$, where $M=(0,L_1)\times(0, L_2)$, with positive constants $L_1, L_2$ and $h$. We complement it with the
following boundary conditions
\begin{eqnarray}
& v, w, T \mbox{ are }\mbox{periodic in }x \mbox{ and }y ,\label{PEBC1'}\\
&(\partial_zv,w)|_{z=-h,0}=0,\label{PEBC2'}\\
&T|_{z=-h}=1,\quad T|_{z=0}=0,\label{PEBC3'}
\end{eqnarray}
and the initial condition
\begin{eqnarray}
&(v,T)|_{t=0}=(v_0, T_0). \label{PEIC'}
\end{eqnarray}
We point out that here we do not need the initial condition on the vertical velocity $w$, as it can be uniquely recovered from the horizontal velocity $v$ from the incompressibility condition.

Replacing $T$ and $p$ by $T+\frac zh$ and $p+\frac{z^2}{2h}$,
respectively, extending the resultant unknowns
$v, w$ and $T$ evenly, oddly and oddly, respectively, with respect
to $z$, and noticing that the periodic subspace
$\mathcal H$, given by
\begin{align*}
  \mathcal H:=&\{(v,w,p,T)|v,w,p\mbox{ and }T\mbox{are spatially periodic in all three variables} \\
  &\mbox{ and are even, odd, even and odd with respect to }z \mbox{ variable},\mbox{ respectively}\},
\end{align*}
is invariant under the dynamical system (\ref{PEFV'}),
one can easily check that system (\ref{PEFV'}), subject to (\ref{PEBC1'})--(\ref{PEIC'}), defined on the domain $\Omega_0=M\times(-h,0)$, is equivalent to the following system, defined on the extended domain $\Omega:=M\times(-h,h)$,
\begin{equation}
  \label{PEFV}
  \left\{
  \begin{array}{l}
    \partial_tv +(v \cdot\nabla_H)v+w\partial_zv-\Delta v+\nabla_H
    p+f_0k\times v=0,\\
    \nabla_H\cdot v+\partial_zw=0,\\
    \partial_zp=T,\\
    \partial_tT+v\cdot\nabla_HT+w\left(\partial_zT+\frac1h\right)
    -\partial_z^2T=0,
  \end{array}
  \right.
\end{equation}
subject to the following new boundary and initial conditions
\begin{eqnarray}
& v, w, p \mbox{ and } T \mbox{ are }\mbox{periodic in }x, y, z,\label{PEBC1}\\
& v\mbox{ and }p \mbox{ are even in }z,\mbox{ and } w\mbox{ and }T\mbox{ are odd in }z,\label{PEBC2}\\
&(v,T)|_{t=0}=(v_0, T_0).   \label{PEIC}
\end{eqnarray}

The boundary condition (\ref{PEBC1}) and the symmetry condition (\ref{PEBC2}) imply $w|_{z=-h}=0$. Therefore, the vertical velocity can be uniquely determined in terms of the horizontal velocity through the incompressibility condition as
\begin{equation}\label{VV}
  w(x,y,z,t)=-\int_{-h}^z\nabla_H\cdot v(x,y,z',t)dz',
\end{equation}
for any $(x,y,z)\in \Omega$ and $t\in(0,\infty)$.

\begin{definition}
  \label{PEFVDEF}
Suppose that $v_0,T_0\in H^2(\Omega)$ are periodic, and are even and odd in $z$, respectively.
Given a positive time $\mathcal T\in(0,\infty)$. A pair $(v,T)$
is called a strong solution to system (\ref{PEFV}), subject to
(\ref{PEBC1})--(\ref{PEIC}), in $\Omega\times(0,\mathcal T)$, if

(i) $v$ and $T$ are spatially periodic, and are even and odd
in $z$, respectively;

(ii) $v$ and $T$ have the regularities
\begin{align*}
&v\in L^\infty(0,\mathcal T; H^2(\Omega))\cap C([0,\mathcal T];
H^1(\Omega))\cap L^2(0,\mathcal T; H^3(\Omega)),\\
&T\in L^\infty(0,\mathcal T; H^2(\Omega))\cap C([0,\mathcal T];
H^1(\Omega)),\quad\partial_zT\in L^2(0,\mathcal T; H^2(\Omega)),\\
&\partial_tv\in L^2(0,\mathcal T; H^1(\Omega)),\quad\partial_tT
\in L^2(0,\mathcal T;H^1(\Omega));
\end{align*}

(iii) $v$ and $T$ satisfy system (\ref{PEFV}) pointwisely, a.e.\,in
$\Omega\times(0,\mathcal T)$, with $w$ given by (\ref{VV}), and fulfill the initial condition (\ref{PEIC}).
\end{definition}

We have the following result on the existence and uniqueness of strong solutions to system (\ref{PEFV}), subject to (\ref{PEBC1})--(\ref{PEIC}), on $\Omega\times(0,\mathcal T)$, for any finite time $\mathcal T$.

\begin{theorem}[cf. \cite{CAOLITITI1}]
  \label{THMPEFV}
Given a positive time $\mathcal T\in(0,\infty)$. Suppose that the periodic functions $v_0,T_0\in H^2(\Omega)$ are even and odd in $z$, respectively, with $\int_{-h}^h\nabla_H\cdot v_0dz=0,$ on $M$.
Then, there is a unique strong solution $(v,T)$ to
system (\ref{PEFV})--(\ref{PEIC}), in $\Omega\times(0,\mathcal T)$. Moreover, the unique strong solution $(v,T)$ is continuously dependent on the initial data.
\end{theorem}

\begin{remark}
Theorem \ref{THMPEFV} generalizes and complements the results in \cite{CAOTITI12}, where the global regularity of strong solutions was
proved, for initial data with some higher regularities than $H^2$.
Note that the local existence and uniqueness of strong solutions were
not established in \cite{CAOTITI12}.
\end{remark}

\subsection{The horizontal diffusivity case with full viscosity}
\label{SECPEFH}
As a counterpart of the result in the previous subsection, it is also interesting to study the following PEs with full viscosity but only horizontal diffusivity
\begin{equation}
  \label{PEFH'}
  \left\{
  \begin{array}{l}
    \partial_tv +(v \cdot\nabla_H)v+w\partial_zv-\Delta v+\nabla_H
    p+f_0k\times v=0,\\
    \nabla_H\cdot v+\partial_zw=0,\\
    \partial_zp=T,\\
    \partial_tT+v\cdot\nabla_HT+w\partial_zT
    -\Delta_HT=0,
  \end{array}
  \right.
\end{equation}
in $\Omega_0=M\times(-h,0)$. We complement it with the boundary and initial conditions (\ref{PEBC1'})--(\ref{PEIC'}). Same to the vertical diffusivity case, by replacing $T$ and $p$ by $T+\frac zh$ and $p+\frac{z^2}{2h}$,
respectively, extending the resultant unknowns
$v, w$ and $T$ evenly, oddly and oddly, respectively, with respect
to $z$, system (\ref{PEFH'}), defined on $\Omega_0$, subject to (\ref{PEBC1'})--(\ref{PEIC'}), is equivalent to the following system:
\begin{equation}
  \label{PEFH}
  \left\{
  \begin{array}{l}
    \partial_tv +(v \cdot\nabla_H)v+w\partial_zv-\Delta v+\nabla_H
    p+f_0k\times v=0,\\
    \nabla_H\cdot v+\partial_zw=0,\\
    \partial_zp=T,\\
    \partial_tT+v\cdot\nabla_HT+w\left(\partial_zT+\frac1h\right)
    -\Delta_HT=0,
  \end{array}
  \right.
\end{equation}
in the extended domain $\Omega=M\times(-h,h)$, subject to the boundary and initial conditions (\ref{PEBC1})--(\ref{PEIC}).

Strong solutions to system (\ref{PEFH}), subject to the boundary and initial conditions (\ref{PEBC1})--(\ref{PEIC}), are defined similarly as Definition \ref{PEFVDEF}. The only difference is that we now ask for the regularity $\nabla_HT\in L^2(0,\mathcal T; H^2(\Omega))$, instead of the regularity $\partial_zT\in L^2(0,\mathcal T; H^2(\Omega))$.

We have the following counterpart of Theorem \ref{THMPEFV}:

\begin{theorem}[cf. \cite{CAOLITITI2}]
  \label{THMPEFH}
Given a positive time $\mathcal T\in(0,\infty)$. Suppose that the periodic functions $v_0,T_0\in H^2(\Omega)$ are even and odd in $z$, respectively, with $\int_{-h}^h\nabla_H\cdot v_0dz=0,$ on $M$.
Then there is a unique strong solution $(v,T)$ to
system (\ref{PEFV}), subject to the boundary and initial conditions (\ref{PEBC1})--(\ref{PEIC}), in $\Omega\times(0,\mathcal T)$. Moreover, the unique strong solution $(v,T)$ is continuously dependent on the initial data.
\end{theorem}


\section{The PEs with horizontal viscosity and partial diffusivity}
\label{SECPEHP}
In the oceanic and atmospheric dynamics, due to the strong horizontal turbulent mixing,
the horizontal viscosity is much stronger than the vertical viscosity.
Therefore, both physically and mathematically, it is important to study the limiting case that the vertical viscosity vanishes.

\subsection{The horizontal diffusivity case with horizontal viscosity}
\label{SECPEHH}
Consider the following PEs with only horizontal viscosity and horizontal diffusivity:
\begin{equation}\label{PEHH}
\left\{
\begin{array}{l}
\partial_tv+(v\cdot\nabla_H)v+w\partial_zv+\nabla_Hp-\Delta_H v+f_0k\times v=0,\\
\partial_zp+T=0,\\
\nabla_H\cdot v+\partial_zw=0,\\
\partial_tT+v\cdot\nabla_HT+w\left(\partial_zT+\frac{1}{h}\right) -\Delta_HT=0,
\end{array}
\right.
\end{equation}
in $\Omega=M\times(-h,h)$, with $M=(0,L_1)\times(0,L_2)$. We complement it with the boundary and initial conditions (\ref{PEBC1})--(\ref{PEIC}).
Let us recall that the boundary conditions (\ref{PEBC1})--(\ref{PEBC2}) are equivalent to the boundary conditions (\ref{PEBC1'})--(\ref{PEBC3'}), by performing suitable transformations to the unknowns.

Strong solutions to (\ref{PEHH}), subject to (\ref{PEBC1})--(\ref{PEIC}), in $\Omega\times(0,\mathcal T)$, are define in the similar way as Definition \ref{PEFVDEF}, by replacing (ii) there with the following:

(ii') $v$ and $T$ have the regularities
\begin{eqnarray*}
&(v,T)\in L^\infty(0,\mathcal T; H^2(\Omega))\cap C([0,\mathcal T];H^1(\Omega)),\\
&(\nabla_Hv,\nabla_HT)\in L^2(0,\mathcal T; H^2(\Omega)), \quad(\partial_tv,\partial_tT)\in L^2(0,\mathcal T; H^1(\Omega)).
\end{eqnarray*}

We have the following theorem on the existence and uniqueness of strong
solutions to system (\ref{PEHH}), subject to (\ref{PEBC1})--(\ref{PEIC}), in $\Omega\times(0,\mathcal T)$, for any finite time $\mathcal T$.

\begin{theorem}[cf. \cite{CAOLITITI3}]
  \label{THMPEHH}
Given a positive time $\mathcal T\in(0,\infty)$. Suppose that the periodic functions $v_0,T_0\in H^2(\Omega)$ are even and odd in $z$, respectively, with $\int_{-h}^h\nabla_H\cdot v_0dz=0,$ on $M$.
Then there is a unique strong solution $(v,T)$ to
system (\ref{PEHH}), subject to the boundary and initial conditions
(\ref{PEBC1})--(\ref{PEIC}), in $\Omega\times(0,\mathcal T)$.
Moreover, the unique strong solution $(v,T)$ is continuously
dependent on the initial data.
\end{theorem}

The key issue of proving Theorem \ref{THMPEHH} is establishing the a priori $H^2$ estimates on the strong solutions, up to any finite time. Our analysis shows that all desired estimates depend on the $L^2$ estimate on $\partial_zv$. Due to the lack of the vertical viscosity in the horizontal momentum equations, one will encounter a factor $\|v\|_\infty^2$ in the energy inequalities, in other words,
the energy inequalities obtained are of the following form:
$$
\frac{d}{dt}f\leq C\|v\|_{\infty}^2f+``\text{other terms}",
$$which forces us somehow
to establish the a priori estimate $\int_0^\mathcal T\|v\|_{\infty}^2dt$. Generally, one can not achieve such kind \emph{a priori} estimate for any finite time. Nevertheless, we have the
following proposition, which states the \emph{a priori} estimate on the growth of the $L^q$-norms, for all $q\in[2,\infty)$, of the velocity $v$.

\begin{proposition}[cf. \cite{CAOLITITI3}]
\label{PROPPEHH}
Let $(v,T)$ be a strong solution to system (\ref{PEHH}), subject to the boundary and initial conditions (\ref{PEBC1})--(\ref{PEIC}), on the interval $(0,\mathcal T)$. Then, for each $q\in[2,\infty)$, the following estimate holds
$$
\sup_{0\leq t\leq\mathcal T}\|v\|_q \leq K_1(\mathcal T)e^{C\|T_0\|_q^2\mathcal T}
(1+\|v_0\|_q)\sqrt q,
$$
where $K_1$ is a continuously increasing function, determined by the initial norms $\|v_0\|_{2}$, $\|T_0\|_{2}$, $\|v_0\|_{4}$ and $\|T_0\|_{4}$.
\end{proposition}

Thanks to the above proposition, we can control the main part of the quantity $\|v\|_{\infty}^2$. To this end, we use a logarithmic Sobolev limiting inequality, generalizing the classical Br\'ezis--Gallout--Wainger inequality \cite{Brezis_Gallouet_1980,Brezis_Wainger_1980} (see also \cite{CAOFARHATTITI}), stated in the next proposition. This logarithmic inequality serves as a bridge between the $L^q$ norms and the $L^\infty$ norm.

\begin{lemma}[cf. \cite{CAOLITITI3}]\label{LOG}
Let $F\in W^{1,p}(\Omega)$, with $p>3$, be a periodic function. Then the following inequality holds
\begin{equation*}
  \|F\|_\infty\leq C_{ {p},\lambda}\max\left\{1,\sup_{r\geq2}\frac{\|F\|_r}{r^\lambda}\right\}
  \log^\lambda
  (\|F\|_{W^{1, {p}}(\Omega)}+e),
\end{equation*}
for any $\lambda>0$.
\end{lemma}

Now, applying Proposition \ref{PROPPEHH} and Lemma \ref{LOG}, we can control $\|v\|_{\infty}^2$ in terms of the logarithm of higher order norms as follows
\begin{equation*}
\|v\|_{\infty}^2\leq K(t)\log (\|(v,T)\|_{H^2}^2+e),
\end{equation*}
for a function $K\in L^1((0,\mathcal T))$. Note that such kind logarithm dependence of $\|v\|_{\infty}^2$ on the higher order norms
does not effect us to obtain the
a priori estimate, up to any finite time. Therefore, one
can successfully achieve the desired a priori $H^2$ estimate, up to any finite time, on strong solution $(v,T)$.

Besides, the following lemma of a system version of the classic Gronwall inequality plays an important role in simplifying the proof of the a priori estimates.

\begin{lemma}[cf. \cite{CAOLITITI3}]
  \label{gronwall}
Let $m(t), K(t),
A_i(t)$ and $B_i(t)$ be nonnegative functions, such that $A_i\geq e,$ are absolutely continuous, for $i=1,\cdots,n, K\in L_{\text{loc}}^1([0,\infty))$.
Given a positive time $\mathcal T$, and suppose that
\begin{align*}
&\frac{d}{dt}A_1(t)+B_1(t)\leq K(t)\left(\log\sum_{i=1}^nA_i(t)\right)A_1(t), \\
&\frac{d}{dt}A_i(t)+B_i(t)\leq K(t)\left(\log\sum_{i=1}^nA_i(t)\right)A_i(t)+\zeta A_{i-1}^\alpha(t)B_{i-1}(t),
\end{align*}
for $i=1,2,\cdots,n$, and for any $t\in(0,\mathcal T)$, where $\alpha\geq1$ and $\zeta\geq1$ are two constants.
Then it holds that
$$
\sum_{i=1}^nA_i(t)+\sum_{i=1}^n\int_0^tB_i(s)ds\leq Q(t),\quad \forall t\in[0,\mathcal T),
$$
where $Q$ is a continuous function on $[0,\infty)$ which is determined by $A_i(0), i=1,\cdots,n,$ and $K$.
\end{lemma}

Some generalizations of Theorem \ref{THMPEHH} has been obtained in \cite{CAOLITITI5}, stated in the following theorem.

\begin{theorem}[cf. \cite{CAOLITITI5}]
\label{THMPEHH'}
Suppose that the periodic functions $v_0,T_0\in H^1(\Omega)$ are
even and odd in $z$, respectively, with $\int_{-h}^h\nabla_H\cdot v_0(x,y,z)dz=0$, for any $(x,y)\in M$.  Then, there is a unique
local strong solution $(v,T)$ to system (\ref{PEHH}), subject to
the boundary and initial conditions (\ref{PEBC1})--(\ref{PEIC}).

Moreover, if we assume in addition that
$$
\partial_zv_0\in L^m(\Omega),\quad (v_0, T_0)\in L^\infty(\Omega),
$$
for some $m\in(2,\infty)$, then the corresponding local strong
solution $(v,T)$ can be extended uniquely to any finite time $\mathcal T\in(0,\infty)$.
\end{theorem}

\begin{remark}
It is unclear if system (\ref{PEHH}), subject to
(\ref{PEBC1})--(\ref{PEIC}), has strong solutions up to any finite
time $\mathcal T\in(0,\infty)$, for any initial data $(v_0, T_0)\in H^1$. Even the local existence and uniqueness part of Theorem \ref{THMPEHH'} does not follow from the standard energy
approach, and in fact, some local in space type energy estimate has been used in the proof of local existence of strong solutions in \cite{CAOLITITI5}.
\end{remark}

\subsection{The vertical diffusivity case with horizontal viscosity}
\label{SECPEHV}
As a counterpart of the previous subsection, we consider in this subsection the following PEs with only horizontal viscosity and vertical diffusivity:
\begin{equation}\label{PEHV}
\left\{
\begin{array}{l}
\partial_tv+(v\cdot\nabla_H)v+w\partial_zv+\nabla_Hp-\Delta_H v+f_0k\times v=0,\\
\partial_zp+T=0,\\
\nabla_H\cdot v+\partial_zw=0,\\
\partial_tT+v\cdot\nabla_HT+w\left(\partial_zT+\frac{1}{h}\right) -\partial_z^2T=0,
\end{array}
\right.
\end{equation}
in $\Omega=M\times(-h,h)$, with $M=(0,L_1)\times(0,L_2)$. Same to the previous subsection, we again complement it with the boundary and initial conditions (\ref{PEBC1})--(\ref{PEIC}).

Existence and uniqueness of strong solutions, up to any finite time $\mathcal T\in(0,\infty)$, to the
PEs with horizontal viscosity and vertical diffusivity, i.e.\,system (\ref{PEHV}), is much more complicated than
that to the PEs with both horizontal viscosity and horizontal diffusivity, i.e.\,system (\ref{PEHH}). In fact, due to the ``mismatching of regularities" between $v$ and $T$ in system (\ref{PEHV}), even the definition of strong solutions is not standard.

In order to give the definition of strong solutions, we introduce the following functions
\begin{eqnarray}\label{defuet}
  u=\partial_zv, \quad\theta=\nabla_H^\perp\cdot v,\quad\eta=\nabla_H\cdot v+\int_{-h}^zTd\xi-\frac{1}{2h}\int_{-h}^h\int_{-h}^zTd\xi dz,
\end{eqnarray}
where $\nabla_H^\perp=(-\partial_y, \partial_x)$.

Strong solutions to system (\ref{PEHV}), subject to (\ref{PEBC1})--(\ref{PEIC}), are defined as follows.

\begin{definition}\label{def1.1}
Given a positive time $\mathcal T$. Let $v_0\in H^2(\Omega)$ and $T_0\in H^1(\Omega)$, with $\int_{-h}^h\nabla_H\cdot v_0(x,y,z) dz=0$ and $\nabla_HT_0\in L^4(\Omega)$, be two periodic functions, such that they are even and odd in $z$, respectively. A pair $(v,T)$ is called a strong solution to system
(\ref{PEHV}), subject to (\ref{PEBC1})--(\ref{PEIC}), in $\Omega\times(0,\mathcal T)$, if

(i) $v$ and $T$ are periodic in $x,y,z$, and they are even and odd in $z$, respectively;

(ii) $v$ and $T$ have the regularities
\begin{eqnarray*}
&&v\in L^\infty(0,\mathcal T; H^2(\Omega))\cap C([0,\mathcal T];H^1(\Omega)),\quad\partial_tv\in L^2(0,\mathcal T; H^1(\Omega)),\\
&&T\in L^\infty(0,\mathcal T; H^1(\Omega))\cap C([0,\mathcal T]; L^2(\Omega)), \quad\partial_t T\in L^2(0,\mathcal T; L^2(\Omega)), \\
&&(\nabla_Hu,\partial_zT)\in L^2(0,\mathcal T; H^1(\Omega)),\quad\nabla_HT\in L^\infty(0,\mathcal T; L^4(\Omega)),\\
&&\eta\in L^2(0,\mathcal T; H^2(\Omega)),\quad \theta\in L^2(0,\mathcal T; H^2(\Omega));
\end{eqnarray*}

(iii) $v$ and $T$ satisfy system (\ref{PEHV}) pointwisely, a.e.\,in $\Omega\times(0,\mathcal T)$, with $w$ given by (\ref{VV}), and fulfill the initial condition (\ref{PEIC}).
\end{definition}

\begin{remark}
\label{remark1.1}
(i) The regularities in Definition \ref{def1.1} seem a little bit nonstandard. This is caused by the ``mismatching" of regularities between the horizontal momentum equation $(\ref{PEHV})_1$ and the temperature equation $(\ref{PEHV})_4$: a term involving the horizontal derivatives of the temperature appears in the horizontal momentum equation, but it is only in the vertical direction that the temperature has dissipation. More precisely, though one can obtain the regularity that $\nabla_H\partial_zv\in L^2(0,\mathcal T; H^1(\Omega))$, which is included in Definition (\ref{def1.1}), we have no reason to ask for the regularity that $\nabla_H^2v\in L^2(0,\mathcal T; H^1(\Omega))$, under the assumption on the initial data in Definition \ref{def1.1}. In fact, recalling the regularity theory for the parabolic system, and checking the  horizontal momentum equation $(\ref{PEHV})_1$, the regularity that $\nabla_H^2v\in L^2(0,\mathcal T; H^1(\Omega))$ appeals to somehow $\nabla_H^2T\in L^2(\Omega\times(0,\mathcal T))$; however, this last requirement can not be fulfilled, because we only have the smoothing effect in the vertical direction for the temperature.

(ii) As stated in (i), one can not expect such regularity that $\nabla_H^2v\in L^2(0,\mathcal T; H^1(\Omega))$. However, with the help of $\eta$ and $\theta$, in (\ref{defuet}), one can expect that some appropriate combinations of $\nabla_Hv$ and $T$ can indeed have the second order spatial derivatives, that is $(\eta,\theta)\in L^2(0,\mathcal T; H^2(\Omega))$, as included in Definition \ref{def1.1}.
\end{remark}

We have the following existence and uniqueness result:

\begin{theorem}[cf.\,\cite{CAOLITITI5}]
\label{THMHV}
Given a positive time $\mathcal T\in(0,\infty)$. Let $v_0\in H^2(\Omega)$ and $T_0\in H^1(\Omega)\cap L^\infty(\Omega)$, with $\int_{-h}^h\nabla_H\cdot v_0(x,y,z) dz=0$ and $\nabla_HT_0\in L^4(\Omega)$, be two periodic functions, such that they are even and odd in $z$, respectively.
Then, system (\ref{PEHV}), subject to (\ref{PEBC1})--(\ref{PEIC}), in $\Omega\times(0,\mathcal T)$, has a unique strong solution $(v,T)$, which is continuously depending on the initial data.

If we assume, in addition, that $T_0\in H^2(\Omega)$, then $(v,T)$ obeys the following additional regularities
\begin{eqnarray*}
&&T\in L^\infty(0,\mathcal T; H^2(\Omega))\cap C(0,\mathcal T; H^1(\Omega)),\quad\partial_tT\in L^2(0,\mathcal T; H^1(\Omega)),\\
&&\nabla_Hv\in L^2(0,\mathcal T; H^2(\Omega)),\quad \partial_zT\in L^2(0,\mathcal T; H^2(\Omega)),
\end{eqnarray*}
for any time $\mathcal T\in(0,\infty)$.
\end{theorem}

\begin{remark}
Generally, if we imposed more regularities on the initial data, then one can expect more  regularities of the strong solutions, and in particular, the strong solution will belong to $C^\infty(\overline\Omega\times[0,\infty))$, as long as the initial datum lies in $C^\infty(\overline\Omega)$. However, one can not expect that the solutions have as high regularities as desired, if the initial data are not accordingly smooth enough.
\end{remark}


The main difficulties for the mathematical analysis of system
(\ref{PEHV}) come from three aspects: one is the strongly nonlinear
term $w\partial_zv$, the next one is the absence of the vertical
viscosity and the horizontal diffusivity, and the last one is the
``mismatching" of regularities between the horizontal momentum
equations and the temperature equation. To deal with the difficulties
caused by the strongly nonlinear term $w\partial_zv$, one can adopt
the ideas that have been frequently used in
\cite{CAOTITI03,CAOTITI07,CAOTITI12,CAOLITITI1,CAOLITITI2,CAOLITITI3},
i.e.\,using the Ladyzhenskaya type inequalities for a class of
integrals in 3D, among which the prototype ones are stated in
the following:

\begin{lemma}[cf. \cite{CAOLITITI3}] \label{lad}
The following inequalities hold true
\begin{align*}
&\int_M\left(\int_{-h}^h|\phi(x,y,z)|dz\right)\left(\int_{-h}^h|\varphi(x,y,z)\psi(x,y,z)|dz\right)dxdy\\
\leq&C\min\left\{\|\phi\|_2^{\frac{1}{2}}\left(\|\phi\|_2^{\frac{1}{2}}
+\|\nabla_H\phi\|_{2}^{\frac{1}{2}
}\right)\|\varphi\|_2\|\psi\|_2^{\frac{1}{2}}\left(
\|\psi\|_2^{\frac{1}{2}}+\|\nabla_H\psi\|_2^{\frac{1}{2}}\right)\right.,\\
&\left.\|\phi\|_2\|\varphi\|_2^{\frac{1}{2}}\left(\|\varphi\|_2^{\frac{1}{2}}+\|\nabla_H\varphi\|_{2}^{\frac{1}{2}}
\right)\|\psi\|_2^{\frac{1}{2}}\left(
\|\psi\|_2^{\frac{1}{2}}+\|\nabla_H\psi\|_2^{\frac{1}{2}}\right)\right\},
\end{align*}
and
\begin{align*}
&\int_M\left(\int_{-h}^h|\phi(x,y,z)|dz\right)\left(\int_{-h}^h|\varphi(x,y,z)
\nabla_H\Psi(x,y,z)|dz\right)dxdy\\
\leq&C\min\left\{\|\phi\|_2^{\frac{1}{2}}\left(\|\phi\|_2^{\frac{1}{2}}
+\|\nabla_H\phi\|_{2}^{\frac{1}{2}
}\right)\|\varphi\|_2\|\Psi\|_\infty^{\frac{1}{2}}\|\nabla_H^2\Psi\|_2^{\frac{1}{2}},\right.\\
&\left.\|\phi\|_2\|\varphi\|_2^{\frac{1}{2}}\left(\|\varphi\|_2^{\frac{1}{2}}+\|\nabla_H\varphi\|_{2}^{\frac{1}{2}}
\right)\|\Psi\|_\infty^{\frac{1}{2}}\|\nabla_H^2\Psi\|_2^{\frac{1}{2}}\right\},
\end{align*}
for every $\phi,\varphi,\psi,\Psi$ such that the right hand sides make sense and are finite. Moreover, if $\phi$ has the form $\phi=\nabla_Hf$, for a function $f$, then by the Poinc\'are inequality, the lower order term $\|\phi\|_2^{\frac{1}{2}}$ in the parentheses can be dropped in the above inequalities, and the same words can be said for $\varphi$ and $\psi$.
\end{lemma}

Concerning the difficulties caused by the absence of the vertical viscosity, we can adopt the idea in the previous subsection. As stated in the previous subsection, the absence of the vertical viscosity forces us to estimate $\|v\|_\infty^2$, which appears as factors on the right hand sides of the energy inequalities. To get this estimate, making use of the growth of the $L^q$ norms of $v$, i.e.\,Proposition \ref{PROPPEHH}, and applying the logarithmic Sobolev limiting embedding
inequality, i.e.\,Lemma \ref{LOG}, we can estimate $\|v\|_{\infty}^2$
as
$$
\|v\|_{\infty}^2\leq C\log (e+\|v\|_{W^{1,4}}).
$$
Recalling the definitions of $u,\eta$ and $\theta$, in (\ref{defuet}), by the elliptic estimates, one can control $\|v\|_{W^{1,4}}$ by the $L^4$-norms of $u, \eta$ and $\theta$. Note that $(u,\eta,\theta)$, defined by (\ref{defuet}), satisfies
\begin{align}
\partial_tu+(v\cdot\nabla_H)u+w\partial_zu-&\Delta_Hu+f_0k\times u\nonumber\\
=&-(u\cdot\nabla_H)v+(\nabla_H\cdot v)u+\nabla_HT, \label{u}\\
\partial_t\eta-\Delta_H\eta+\nabla_H\cdot[(v\cdot\nabla_H)v  &+w\partial_zv+f_0k\times v]\nonumber\\
=&-\partial_zT+wT+\int_{-h}^z\nabla_H\cdot(vT)d\xi+f(x,y,t), \label{eta}\\
\partial_t\theta-\Delta_H\theta=&-\nabla_H^\perp\cdot[(v\cdot\nabla_H)v  +w\partial_zv+f_0k\times v], \label{theta}
\end{align}
where $f=f(x,y,t)$ is given by
\begin{align*}
  f=&\frac{1}{2h}\int_{-h}^h\left(\int_{-h}^z\nabla_H\cdot(vT)d\xi+wT+\nabla_H\cdot \big(\nabla_H\cdot(v\otimes v)+f_0k\times v\big)\right)\label{ff}
   dz.
\end{align*}
By performing the $L^2$ energy estimate to the above system for $(u,\eta,\theta)$, and with the help of the
control on $\|v\|_\infty^2$ stated in the above, we can obtain the a
priori ``\,lower order estimates\,". Actually, we have the following:

\begin{proposition}[cf. \cite{CAOLITITI4}]\label{PROP1PEHV}
Let $(v,T)$ be a strong solution to system (\ref{PEHV}), subject to (\ref{PEBC1})--(\ref{PEIC}), with initial data $(v_0, T_0)$ satisfying the assumptions in Theorem \ref{THMHV}. Let $\mathcal T^*$ be the maximal existence time of $(v,T)$. Then, we have the estimate
\begin{align*}
  \sup_{0\leq s\leq t}(\|\theta\|_2^2+\|\eta\|_2^2+\|u\|_4^4)+\int_0^t(\|\nabla_H\theta\|_2^2 +\|\nabla_H\eta\|_2^2+\|\nabla_Hu\|_2^2)ds\leq M_1,
\end{align*}
for any $t\in(0,\mathcal T^*)$, where $M_1$ is a positive constant.
\end{proposition}


To obtain
the strong solutions up to any finite time,
we still need to obtain ``\,higher
order estimates\,". When doing the energy inequalities for the
horizontal derivatives of $T$, caused by the absence of the horizontal
diffusivity in the temperature equation, one need appeal to somehow
$L^\infty$ estimate on $\nabla_Hv$ to deal with the hardest term
$\int_\Omega|\nabla_Hv||\nabla_HT|^4dxdydz.$ To deal with such kind term, we decompose $v$ as
\begin{equation*}\label{zeta}
v(x,y,z,t)=\zeta(x,y,z,t)-\varpi(x,y,z,t),
\end{equation*}
where $\varpi(\cdot,z,t)$ is the unique solution to the elliptic system
\begin{equation*}
\left\{
\begin{array}{l}
\nabla_H\cdot\varpi(x,y,z,t)=\Phi(x,y,z,t)-\frac{1}{|M|}\int_M\Phi(x,y,z,t) dxdy,\quad\mbox{in }\Omega,\\
\nabla_H^\perp\cdot\varpi(x,y,z,t)=0,\quad\mbox{in }\Omega,\qquad\int_M\varpi(x,y,z,t) dxdy=0,
\end{array}
\right.\label{beta}
\end{equation*}
where $\Phi$ is the function given by
\begin{equation*}\label{Phi}
\Phi(x,y,z,t)=\int_{-h}^zT(x,y,\xi,t)d\xi -\frac{1}{2h}\int_{-h}^h\int_{-h}^zT(x,y,\xi,t)d\xi dz.
\end{equation*}
Roughly speaking, $\varpi$ captures the temperature-dependent part of the velocity, while $\zeta$ captures the temperature-independent part of the velocity. Recalling the definitions of $\eta$ and $\theta$, one can easily check that
\begin{equation*}\label{zeeithta}
 \nabla_H\cdot\zeta= \eta-\frac{1}{|M|}\int_M\Phi dxdy,\quad\nabla_H^\perp\cdot\zeta=\theta.
\end{equation*}

With the help of the above decomposition on $v$, performing higher
order energy estimates to system (\ref{u})--(\ref{theta}), and
using the
classical logarithmic Sobolev limiting inequalities of the
Br\'ezis-Gallouet-Wainger \cite{Brezis_Gallouet_1980,Brezis_Wainger_1980} and Beale-Kato-Majda \cite{BKM} types to $\zeta$ and $\varpi$, respectively, one can deal with the the hardest term
$\int_\Omega|\nabla_Hv||\nabla_HT|^4dxdydz$, and thus get the higher order energy inequality. We have the following proposition:

\begin{proposition}[cf. \cite{CAOLITITI4}]
 \label{PROP2PEHV}
Let $(v,T)$ be a strong solution to system (\ref{PEHV}), subject to (\ref{PEBC1})--(\ref{PEIC}), with initial data $(v_0, T_0)$ satisfying the assumptions in Theorem \ref{THMHV}. Let $\mathcal T^*$ be the maximal existence time of $(v,T)$. Then, we have the estimate
\begin{align*}
  &\sup_{0\leq s\leq t}(\|\nabla_H\eta\|_2^2+\|\nabla_H\theta\|_2^2+\|\nabla u\|_2^2+\|\nabla T\|_2^2+{\|\nabla_HT\|_4^4})\\
  &+\int_0^t(\|\Delta_H\eta\|_2^2+\|\Delta_H\theta\|_2^2 +\|\nabla_H\nabla u\|_2^2+\|\nabla\partial_z T\|_2^2)ds\leq M_2,
\end{align*}
for any $t\in(0,\mathcal T^*)$, where $M_2$ is a positive constant.
\end{proposition}

With the help of Proposition \ref{PROP1PEHV} and Proposition \ref{PROP2PEHV}, it is then standard to show the existence and uniqueness of strong solutions to system (\ref{PEHV}), subject to (\ref{PEBC1})--(\ref{PEIC}), up to any finite time $\mathcal T\in(0,\infty)$.


\section{Rigorous justification of the hydrostatic balance}
\label{SECPEHYD}
Recall that the hydrostatic balance, i.e.\,(\ref{HB}), is a fundamental equation in the PEs, which is derived by taking the formal small aspect ratio limit to the Boussinesq equations (\ref{BE}), as stated in section \ref{SECPEINT}. This section is devoted to the rigorous justification of this small aspect ratio limit. For simplicity, we
ignore the temperature equation and the rotating term $f_0k\times v$; however, the same results hold for the general case, and they can be
proved similarly.

Consider the anisotropic Navier-Stokes equations
$$
\partial_t\mathcal U+(\mathcal U\cdot\nabla)\mathcal U-\mu\Delta_H\mathcal U-\nu\partial_z^2\mathcal U+\nabla p=0,
$$
in the $\varepsilon$-dependent domain
$\Omega_\varepsilon:=M\times(-\varepsilon,\varepsilon)$, where
$M=(0,L_1)\times(0,L_2)$.
Here $\mathcal U=(v, w)$, with $v=(v^1, v^2)$. As in section \ref{SECPEINT}, we suppose that $\mu=1$ and $\nu=\varepsilon^2$.
Introduce the new scaled unknowns
\begin{eqnarray*}
  &&u_\varepsilon=(v_\varepsilon, w_\varepsilon),\quad
  v_\varepsilon(x,y,z,t)=v(x,y,\varepsilon z,t),\\
  &&w_\varepsilon(x,y,z,t)=\frac{1}{\varepsilon}w(x,y,\varepsilon z,t),\quad
  p_\varepsilon(x,y,z,t)=p(x,y,\varepsilon z,t),
\end{eqnarray*}
for any $(x,y,z)\in \Omega:=M\times(-1,1)$, and for any $t\in(0,\infty)$.
Then, $u_\varepsilon=(v_\varepsilon, w_\varepsilon)$ and $p_\varepsilon$ satisfy
the following scaled Navier-Stokes equations (SNS)
\begin{equation}\label{SNS}
(SNS)~~\left\{
\begin{array}{l}
  \partial_t v_\varepsilon+(v_\varepsilon\cdot\nabla_H)v_\varepsilon +w_\varepsilon\partial_zv_\varepsilon-\Delta v_\varepsilon+\nabla_Hp_\varepsilon=0,\\
  \nabla_H\cdot v_\varepsilon+\partial_zw_\varepsilon=0,\\
  \varepsilon^2(\partial_tw_\varepsilon+v_\varepsilon\cdot\nabla_H w_\varepsilon+w_\varepsilon\partial_zw_\varepsilon-\Delta w_\varepsilon)+\partial_zp_\varepsilon=0,
\end{array}
\right.
\end{equation}
defined in the fixed domain $\Omega$.

We complement (\ref{SNS}) with the boundary and initial conditions
\begin{eqnarray}
  &v_\varepsilon, w_\varepsilon \mbox{ and }p_\varepsilon\mbox{ are periodic in }x,y,z,\label{bc}\\
&v_\varepsilon, w_\varepsilon\mbox{ and }p_\varepsilon \mbox{ are even, odd and odd in }z,\mbox{ respectively},\label{sc}\\
  &(v_\varepsilon, w_\varepsilon)|_{t=0}=(v_0, w_0). \label{ic}
\end{eqnarray}
Note that (\ref{sc}) is a symmetry condition, which is preserved by (\ref{SNS}), in other words,
it is automatically satisfied for all time, as long as it is
satisfied initially.

%
Formally, by taking the limit $\varepsilon\rightarrow0$ in (SNS), one obtains the following PEs
\begin{equation}\label{PE}
(PEs)~\left\{
\begin{array}{l}
  \partial_tv+(v\cdot\nabla_H)v+w\partial_zv-\Delta v+\nabla_Hp=0,\\
  \nabla_H\cdot v+\partial_zw=0,\\
  \partial_zp=0.
\end{array}
\right.
\end{equation}


Recall that the solutions under consideration satisfy the symmetry condition (\ref{sc}),
so does the initial datum $u_0=(v_0,w_0)$. Since $w_0$ is odd in $z$, one has $w_0|_{z=0}=0$. Thus,
it follows from the incompressibility condition that
\begin{equation}\label{ne00}
w_0(x,y,z)=-\int_0^z\nabla_H\cdot v_0(x,y,z')dz',
\end{equation}
for any $(x,y)\in M$ and $z\in(-1,1)$. Due to this fact, we only need to specify the horizontal components $v_0$, while the vertical component
$w_0$ is uniquely determined in terms of $v_0$ through (\ref{ne00}).

In case that the initial data $v_0\in H^1(\Omega)$, one can not generally
expect that $w_0$, determined by (\ref{ne00}), belongs to $H^1(\Omega)$. Instead, one should consider
$u_0=(v_0,w_0)$ as
an element in $L^2(\Omega)$, and thus can only obtain a global weak solution $(v_\varepsilon, w_\varepsilon)$ to (\ref{SNS}), subject
to (\ref{bc})--(\ref{sc}). For this case, we have the following theorem concerning the strong convergence:

\begin{theorem}[cf. \cite{LITITIHYDRO}]\label{thm0}
Given a periodic function $v_0\in H^1(\Omega)$, such that it is even in $z$, and
$$
\nabla_H\cdot\left(\int_{-1}^1v(x,y,z)dz\right)=0,\quad\int_\Omega v_0(x,y,z) dxdydz=0.
$$
Let $(v_\varepsilon,
w_\varepsilon)$ and $(v,w)$, respectively, be an arbitrary weak
solution to (SNS) and the unique global strong solution to (PEs), subject to (\ref{bc})--(\ref{sc}). Denote
$$
(V_\varepsilon, W_\varepsilon)=(v_\varepsilon-v, w_\varepsilon-w).
$$

Then, we have the a priori estimate
\begin{align*}
  \sup_{0\leq t<\infty}\|(V_\varepsilon,\varepsilon W_\varepsilon)\|_2^2) +\int_0^{\infty}\|\nabla( V_\varepsilon,\varepsilon W_\varepsilon) \|_2^2dt
  \leq C\varepsilon^2(\|v_0\|_2^2+\varepsilon^2\|w_0\|_2^2+1)^2,
\end{align*}
for any $\varepsilon\in(0,\infty)$, where $C$ is a positive constant depending only on $\|v_0\|_{H^1}$, $L_1$, and $L_2$. As a consequence, we have the following strong convergences
\begin{eqnarray*}
  &(v_\varepsilon,\varepsilon w_\varepsilon)\rightarrow (v,0),\mbox{ in }L^\infty(0,\infty; L^2(\Omega)),\\
  &(\nabla v_\varepsilon, \varepsilon\nabla w_\varepsilon,w_\varepsilon)\rightarrow(\nabla v,0,w),\mbox{ in }L^2(0,\infty;L^2(\Omega)),
\end{eqnarray*}
and the convergence rate is of the order $O(\varepsilon)$.
\end{theorem}

If, in addition, we suppose that $v_0\in H^2(\Omega)$, then $u_0=(v_0,w_0)\in H^1(\Omega)$,
with $w_0$ given by (\ref{ne00}). For this case, we have the following theorem concerning the strong convergence, in which the convergence
is stronger than
that in Theorem \ref{thm0}:

\begin{theorem}[cf. \cite{LITITIHYDRO}]
  \label{thm}
Given a periodic function $v_0\in H^2(\Omega)$, such that it is even in $z$, and
$$
\nabla_H\cdot\left(\int_{-1}^1v(x,y,z)dz\right)=0,\quad\int_\Omega v_0(x,y,z) dxdydz=0.
$$
Let $(v_\varepsilon,
w_\varepsilon)$ and $(v,w)$, respectively, be the unique local strong solution to (SNS) and the unique global strong
solution to (PEs), subject to (\ref{bc})--(\ref{sc}). Denote
$$
(V_\varepsilon, W_\varepsilon)=(v_\varepsilon-v, w_\varepsilon-w).
$$

Then, there is a positive constant $\varepsilon_0$ depending only on the initial norm $\|v_0\|_{H^2}$, $L_1$ and $L_2$, such that, for any
$\varepsilon\in(0,\varepsilon_0)$, the strong solution $(v_\varepsilon, w_\varepsilon)$ exists globally in time, and the following estimate holds
$$
\sup_{0\leq t<\infty}\|(V_\varepsilon,\varepsilon W_\varepsilon)\|_{H^1}^2+
\int_0^\infty\|\nabla(V_\varepsilon,\varepsilon W_\varepsilon) \|_{H^1}^2dt\leq C\varepsilon^2,
$$
for a constant $C$ depending only on $\|v_0\|_{H^2}$, $L_1$ and $L_2$.
As a consequence, the following strong convergences hold
\begin{eqnarray*}
  &(v_\varepsilon,\varepsilon w_\varepsilon)\rightarrow (v,0),\mbox{ in }L^\infty(0,\infty; H^1(\Omega)),\\
  &(\nabla v_\varepsilon, \varepsilon\nabla w_\varepsilon,w_\varepsilon)\rightarrow(\nabla v,0,w),\mbox{ in }L^2(0,\infty;H^1(\Omega)),\\
  &w_\varepsilon\rightarrow w,\mbox{ in }L^\infty(0,\infty;L^2(\Omega)),
\end{eqnarray*}
and the convergence rate is of the order $O(\varepsilon)$.
\end{theorem}

\begin{remark}
(i) Weak convergence from (SNS) to (PEs), as the aspect ratio goes to zero, was established in \cite{AZGU}. However, since only the weak convergence
was obtained in \cite{AZGU}, no convergence rate was provided there.
While Theorems \ref{thm0} and \ref{thm} show that the strong convergence from (SNS) to
(PEs) is global and uniform in time, and the convergence rate is of the
same order to the aspect ratio $\varepsilon$.

(ii) The assumption $\int_\Omega v_0dxdydz=0$
is imposed only for the simplicity of the proof, and the same result
still holds for the general case. One can follow the proof presented
in \cite{LITITIHYDRO}, and do the the a priori
estimates on $v_\varepsilon-v_{0\Omega}$ and $v-v_{0\Omega}$, in stead of on $v_\varepsilon$ and $v$ themselves,
where $v_{0\Omega}=\int_\Omega v_0dxdydz$.


(iii) Generally, if $v_0\in H^k$, with $k\geq2$, then
  $$
  \sup_{0\leq t<\infty}\|(V_\varepsilon,\varepsilon W_\varepsilon)\|_{H^{k-1}}^2+
  \int_0^\infty\|\nabla(V_\varepsilon,\varepsilon W_\varepsilon) \|_{H^{k-1}}^2dt\leq C\varepsilon^2,
  $$
for a positive constant $C$ depending only on $\|v_0\|_{H^k}$, $L_1$ and $L_2$. This can
be done by carrying out higher energy estimates to the difference system between (SNS) and (PEs).

(iv) Observing the smoothing effects of the (SNS) and (PEs) to the strong solutions, one can also
show, in Theorem \ref{thm} (but not in Theorem \ref{thm0}), the strong convergence in higher spaces,
away from the initial time,
in particular, $(v_\varepsilon, w_\varepsilon)\rightarrow(v,w)$, in $C^k(\overline\Omega
\times(T,\infty))$, for any given positive time $T$ and nonnegative integer $k$.
\end{remark}

The proofs of Theorem \ref{thm0} and Theorem \ref{thm} consist of two main ingredients:
the a priori estimates
on the global strong solution $(v,w)$ to
(PEs), and the a priori estimates on the difference $U_\varepsilon=(V_\varepsilon, W_\varepsilon):
=(v_\varepsilon, w_\varepsilon)-(v, w)$.
The desired a priori estimates for strong solution $(v,w)$ to (PEs), subject to (\ref{bc})--(\ref{sc}) is stated in the following proposition.

\begin{proposition}[cf.\,\cite{LITITIHYDRO}]
  \label{cor}
Suppose that $v_0\in H^m$, with $m=1$ or $m=2$. Let $(v,w)$ be the unique global strong solution to (PEs), subject to (\ref{bc})--(\ref{sc}). Then, we have the following:

(i) If $v_0\in H^1(\Omega)$, then we have the estimate
$$
\sup_{0\leq t<\infty}\|v\|_{H^1}^2+\int_0^\infty(\|\nabla v\|_{H^1}^2+
\|\partial_tv\|_2^2)dt\leq C(\|v_0\|_{H^1}, L_1, L_2);
$$

(ii) If $v_0\in H^2(\Omega)$, then we have the estimate
$$
\sup_{0\leq t<\infty}\|v\|_{H^2}^2+\int_0^\infty(\|\nabla v\|_{H^2}^2+
\|\partial_tv\|_{H^1}^2)dt\leq C(\|v_0\|_{H^2}, L_1, L_2).
$$
\end{proposition}

The treatments on the estimates of the difference function $U_\varepsilon$
are different in the proofs of
Theorem \ref{thm0} and Theorem \ref{thm}. For the case of Theorem \ref{thm0}, since $(v_\varepsilon, w_\varepsilon)$
is only a weak solution, one can only perform the energy estimates in the framework of
the weak solutions. We adopt the idea, which was introduced
by Serrin \cite{SERRIN} to prove the weak-strong uniqueness
of the Navier-Stokes equations; however, the difference for our case is that,
the ``strong solutions" are now those for (PEs), but the
``weak solutions" are those for (SNS), or intuitively, we are somehow doing the weak-strong uniqueness
between two different systems.
Precisely, we will: (i) use
$(v,w)$ as the testing functions for (SNS);
(ii) test (PEs) by $v_\varepsilon$; (iii) perform the basic energy identity of (PEs); (iv)
use the energy inequality for (SNS). Then, we get the desired a priori
estimates for $U_\varepsilon$, stated in Theorem \ref{thm0}.

For the case of Theorem \ref{thm}, since the solutions considered are strong ones, one can
get the desired global in time estimates on $U_\varepsilon$ by using the energy approach to
the system for $U_\varepsilon$, which reads as
 \begin{eqnarray*}
  &\partial_tV_\varepsilon+(U_\varepsilon\cdot\nabla)V_\varepsilon -\Delta V_\varepsilon+\nabla_HP_\varepsilon+( u\cdot\nabla)V_\varepsilon+(U_\varepsilon\cdot\nabla) v=0,\\
  &\nabla_H\cdot V_\varepsilon+\partial_zW_\varepsilon=0,\\
  &\varepsilon^2(\partial_tW_\varepsilon+U_\varepsilon\cdot\nabla W_\varepsilon-\Delta W_\varepsilon+U_\varepsilon\cdot\nabla w+ u\cdot\nabla W_\varepsilon)+\partial_zP_\varepsilon\\
  &=-\varepsilon^2(\partial_tw+u\cdot\nabla w-\Delta  w).
  \end{eqnarray*}
However, some arguments and justifications are required. First, since the initial value of $(V_\varepsilon, W_\varepsilon)$
vanishes,
and there is a small coefficient $\varepsilon^2$ in the front of
the external forcing terms in the above system,
one can perform the energy
approach and adopting the smallness argument to
get the desired a priori estimate on $(V_\varepsilon,
W_\varepsilon)$. In particular, we have the following proposition:

\begin{proposition}[cf. \cite{LITITIHYDRO}]
  \label{firstdiff}
There exists a positive constant $\delta_0$, depending only on $L_1$ and $L_2$, such that, the following estimate holds
  \begin{align*}
    &\sup_{0\leq s\leq t}(\|\nabla V_\varepsilon\|_2^2+\varepsilon^2\|\nabla W_\varepsilon\|_2^2)+\int_0^t(\|\Delta V_\varepsilon\|_2^2+\varepsilon^2\|\Delta W_\varepsilon \|_2^2) ds\\
    \leq&C\varepsilon^2e^{C(1+\varepsilon^4)\int_0^t\|\Delta  v\|_2^2\|\nabla\Delta  v\|_2^2ds}\int_0^t(1+\|\Delta  v\|_2^2) (\|\nabla\partial_t  v\|_2^2+\|\nabla\Delta  v\|_2^2)ds,
  \end{align*}
  for any $t\in[0,T_\varepsilon^*)$, as long as
  $$
  \sup_{0\leq s\leq t}(\|\nabla V_\varepsilon\|_2^2+\varepsilon^2\|\nabla W_\varepsilon\|_2^2)\leq\delta_0^2,
  $$
  where $T_\varepsilon*$ is the maximal existence time of $(v_\varepsilon, w_\varepsilon)$, and $C$ is a positive constant depending only on $L_1$ and $L_2$.
\end{proposition}

Second, one has to pay attention that, when
performing the energy estimates to the horizontal momentum equations for $V_\varepsilon$, no more information
about $W_\varepsilon$ can be used other than that comes from the
incompressibility condition, in other words, the term
$W_\varepsilon\partial_zV_\varepsilon$, in the horizontal momentum
equations for $V_\varepsilon$, can be only dealt with by expressing $W_\varepsilon$ as
$$
W_\varepsilon(x,y,z,t)=-\int_{0}^z\nabla_H\cdot V_\varepsilon(x,y,z',t)dz'.
$$
This is because the explicit dynamical information of $W_\varepsilon$, which comes from the vertical momentum equation,
is always
tied up with the parameter $\varepsilon$, which will finally go to
zero, in other words, the vertical momentum equation for $W_\varepsilon$ provides no $\varepsilon$-independent dynamical information of $W_\varepsilon$.
Keeping these in mind, we can obtain the desired a priori estimate stated in Theorem \ref{thm}.

\section{The PEs with a certain class of discontinuous initial data}
\label{SECPEUNIDIS}

In this section, we present some result on the conditional uniqueness of weak solutions to the PEs, which in particular implies the global existence and uniqueness of weak solutions to the PEs, with a certain class of discontinuous initial data.

For simplicity, we consider the following version of the PEs:
\begin{equation}\label{5PES}
\left\{
\begin{array}{l}
\partial_tv+(v\cdot\nabla_H)v+w\partial_zv+\nabla_Hp(x,y,t)-\Delta v+f_0k\times
v=0,\\
\nabla_H\cdot v+\partial_zw=0,
\end{array}
\right.
\end{equation}
in the
spatial domain $\Omega=M\times(-h, 0)$, with $M=(0,L_1)\times(0,L_2)$. We complement (\ref{5PES}) with the following boundary and initial conditions
\begin{eqnarray}
  &v, w\text{ and }p\text{ are periodic in }x,y\text{ and }z,\label{5BC1}\\
  &v\text{ and }w\text{ are even and odd in }z,\text{ respectively},\label{5BC2}\\
  &v|_{t=0}=v_0. \label{5IC}
\end{eqnarray}

Note that, as before, the vertical velocity $w$ can be uniquely expressed in terms of the horizontal velocity $v$, through the incompressibility condition $(\ref{5PES})_2$ and the symmetry condition (\ref{5BC2}), namely
$$
w(x,y,z,t)=-\nabla_H\cdot\left(\int_0^z v(x,y,\xi,t)d\xi\right).
$$
As a result, system (\ref{5PES}), subject to (\ref{5BC1})--(\ref{5IC}), is equivalent to the following
\begin{equation}
\label{5PESM}
\left\{
\begin{array}{l}
\partial_tv+(v\cdot\nabla_H)v+w\partial_zv+\nabla_Hp(\textbf{x}^H,t)-\Delta v+f_0k\times
v=0,\\
\nabla_H\cdot\left(\int_{-h}^h v(x,y,z,t)dz\right)=0,\\
w(x,y,z,t)=-\nabla_H\cdot\left(\int_0^z v(x,y,\xi,t)d\xi\right),
\end{array}
\right.
\end{equation}
in the domain $\Omega=M\times(-h,h)$, subject to the boundary and initial conditions
\begin{eqnarray}
  &v\mbox{ is periodic in }x,y \text{ and }z, \mbox{ and is even in }z,\label{bc1}\\
  &v|_{t=0}=v_0. \label{5ic}
\end{eqnarray}

Define two spaces $X$ and $\mathcal H$ as
\begin{equation*}\label{X}
X=\{v\in L^6(\Omega)|v\text{ is periodic in }z, \mbox{ and }\partial_zv\in L^2(\Omega)\},
\end{equation*}
and
\begin{align*}\label{H}
\mathcal H=\bigg\{& v\in L^2(\Omega)\bigg|v\mbox{ is periodic in }(x,y,z), \mbox{ even in }z, \\
&\mbox{ and satisfies } \nabla_H\cdot \left(\int_{-h}^hv(x,y,z)dz\right)=0\bigg\}.
\end{align*}
%
%
%
%
%

For any initial data $v_0\in\mathcal H$, following the arguments in
\cite{LTW92A,LTW92B,LTW95}, there is a global weak solution to system
(\ref{5PESM})--(\ref{5ic}); however, as we mentioned before, in section \ref{SECPEINT}, it is still
an open question to show the uniqueness of weak solutions to the
PEs. Nevertheless, we can prove the following theorem
on the conditional uniqueness of weak solutions to the PEs.

\begin{theorem}[cf.\,\cite{LITITIUNIDIS}]
  \label{thmwsu}
  Let $v$ be a global weak solution to system (\ref{5PESM})--(\ref{5ic}). Suppose that there is a positive time $T_v$, such that
  \begin{eqnarray*}
  v(\textbf{x},t)=\bar v(\textbf{x},t)+V(\textbf{x},t),\quad \textbf{x}\in\Omega, t\in(0,T_v),\\
  \partial_z\bar v\in L^\infty(0,T_v; L^2(\Omega))\cap L^2(0,T_v; H^1_\text{per}(\Omega)),\quad V\in L^\infty(\Omega\times(0, T_v)).
  \end{eqnarray*}

  Then, there is a positive constant $\varepsilon_0$ depending only on $h$, such that $v$ is the unique global weak solution to system (\ref{5PESM})--(\ref{5ic}), with the same initial data as $v$, provided
  $$
  \sup_{0<t<T_v'}\|V\|_\infty\leq\varepsilon_0,
  $$
  for some $T_v'\in(0,T_v)$.
\end{theorem}

To prove Theorem \ref{thmwsu}, we first show that any weak solution
to the primitive equations is smooth in the space-time domain $\bar\Omega\times(0,\infty)$, that is we have the following proposition.

\begin{proposition}[cf. \cite{LITITIUNIDIS}]
\label{corinterior}
Let $v$ be a weak solution to system (\ref{5PESM})--(\ref{5ic}). Then, $v$ is smooth away from the initial time, i.e.~$v\in C^\infty(\bar\Omega\times(0,\infty))$.
\end{proposition}

The above proposition can be intuitively seen by noticing that a
weak solution has
$H^1$ regularity, at almost any time immediately after the initial time, and
recalling that, for any $H^1$ initial data, there is a unique global
strong solution to the PEs; however, in order to
rigorously prove this fact, we need the the weak-strong uniqueness
result for the PEs, see the next proposition,
below, where we adopt the idea of Serrin \cite{SERRIN}.

\begin{proposition}[cf. \cite{LITITIUNIDIS}]
  \label{apppropwsu}
  Given the initial data $v_0\in H^1_{\text{per}}(\Omega)\cap
  \mathcal H$. Let $v_{\text{s}}$ and $v_{\text{w}}$ be the unique
  global strong solution and an arbitrary global weak solution,
  respectively, to system (\ref{5PESM})--(\ref{5ic}), with the same
  initial data $v_0$. Then, we have $v_{\text{s}}\equiv v_{\text{w}}$.
\end{proposition}

Since weak
solutions are smooth immediately after the initial time, one can perform the
energy estimates to the difference system between two weak solutions,
on any finite time interval away from the initial time. Next, focusing on the short time interval $(0,T_v)$, and using
the decomposition stated in the theorem, and handling the nonlinear
terms involving each parts in their own ways, we then achieve the
uniqueness.

It should be pointed out that one can not expect that all weak
solutions to the PEs, with general initial data in
$\mathcal H$, have the decomposition stated in the above theorem.
In fact, by the weakly lower semi-continuity of the norms, in order
to have such decomposition, it is necessary to require that the
initial data $v_0$ has the decomposition $v_0=\bar v_0+V_0$, with
$\partial_z\bar v_0\in L^2$ and $V_0\in L^\infty$. Observing this, it is
necessary to state the following theorem on the global existence and
uniqueness of weak solutions to the PEs with such kind of
initial data.

\begin{theorem}[cf. \cite{LITITIUNIDIS}]\label{thmain}
Suppose that the initial data $v_0=\bar v_0+V_0$, with $\bar v_0\in
X\cap\mathcal H$ and $V_0\in L^\infty(\Omega)\cap\mathcal H$. Then the
following hold:

(i) There is a global weak solution $v$ to system (\ref{5PESM})--(\ref{5ic}),
such that
\begin{eqnarray*}
v(\textbf{x},t)=\bar v(\textbf{x},t)+V(\textbf{x},t), \quad\textbf{x}\in\Omega, t\in(0,\infty),\\
\partial_z\bar v\in L_{loc}^\infty([0,\infty);L^2(\Omega))\cap L_{loc}^2([0,\infty); H^1_\text{per}(\Omega)),\\
\sup_{0<s<t}\|V\|_\infty(s)\leq \mu(t)\|V_0\|_\infty,\quad t\in(0,\infty),
\end{eqnarray*}
where
$$
\mu(t)=C_0(1+\|v_0\|_4)^{40}(t+1)^2\exp\{C_0e^{2t}(t+1)(1+\|v_0\|_4)^4\},
$$
for some positive constant $C_0$ depending only on $h$;

(ii) Let $\varepsilon_0$ be the positive constant in Theorem \ref{thmwsu}, then the above weak solution is unique, provided
$\mu(0)\|V_0\|_\infty\leq\frac{\varepsilon_0}{2}$.
\end{theorem}

The uniqueness part of Theorem \ref{thmain} is a direct consequence
of the estimate in (i) and Theorem \ref{thmwsu}.
So the key ingredient of the proof of
Theorem \ref{thmain} is to find the required decomposition.
To this end, we regularize the initial data $v_0$, and solve system
(\ref{5PESM})--(\ref{5ic}), with the regularized initial data
$v_{0\varepsilon}$, obtaining a sequence of solutions $v_\varepsilon$,
which approach $v$, as $\varepsilon\rightarrow0$. Due to the lower
semi-continuity of norms, it suffices to find the corresponding
decomposition of $v_\varepsilon$. We decompose $v_{\varepsilon}$ as
$$
v_\varepsilon=\bar v_\varepsilon+V_\varepsilon,
$$
where $V_\varepsilon$ is the unique solution to the following linear system
\begin{eqnarray}
  &\partial_tV_\varepsilon+(v_\varepsilon\cdot\nabla_H)V_\varepsilon+w_\varepsilon\partial_zV_\varepsilon +\nabla_HP_\varepsilon(\textbf{x}^H ,t)-\Delta V_\varepsilon+f_0k\times V_\varepsilon=0,\label{3.5}\\
  &\int_{-h}^h\nabla_H\cdot V_\varepsilon(\textbf{x}^H,z,t)dz=0,\label{3.6}
\end{eqnarray}
subject to the periodic boundary condition, with initial data $V_{0\varepsilon}$.

The desired estimates on $V_\varepsilon$ and $\bar v_\varepsilon$ are stated in the next proposition.

\begin{proposition}[cf. \cite{LITITIUNIDIS}]
We have the following estimates
\begin{eqnarray*}
\sup_{0\leq s\leq t}\|V_\varepsilon\|_\infty(s)\leq K_3(t)\|V_0\|_\infty,\\
\sup_{0\leq s\leq t}(\|\bar v_\varepsilon\|_2^2+\|\partial_z\bar v_\varepsilon\|_2^2) +\int_0^t(\|\nabla \bar v_\varepsilon\|_2^2+\|\nabla\partial_z\bar v_\varepsilon\|_2^2)d\tau\leq K_4(t),
\end{eqnarray*}
where $K_3$ and $K_4$ are continuous function on $[0,\infty)$, determined only by $h$, and the initial norms $\|v_0\|_6, \|\bar v_0\|_2$ and $\|V_0\|_\infty$.
 \end{proposition}

As an application of Theorem \ref{thmain}, we have the following result, which generalizes the results in \cite{BGMR03,TACHIM,PTZ09,KPRZ}.

\begin{corollary}[cf.\,\cite{LITITIUNIDIS}]\label{thm1.2}
For any $v_0\in X\cap\mathcal H$ (or $v_0\in C_\text{per}(\bar\Omega)\cap\mathcal H$), there is a constant $\sigma_0$, depending only on $h$ and the upper bound of $\|v_0\|_4$, such that for any $ \mathscr V_0=v_0+V_0$, with $V_0\in L^\infty(\Omega)\cap\mathcal H$ and $\|V_0\|_\infty\leq\sigma_0$, system (\ref{5PESM})--(\ref{5ic}) with initial data $\mathscr V_0$ has a unique weak solution, which has the regularities stated in Theorem \ref{thmain}.
\end{corollary}

\begin{example}[cf.\,\cite{LITITIUNIDIS}]
Given a constant vector $a=(a^1, a^2)$, and two positive numbers
$\delta$ and $\eta$, with $\eta\in(0,h)$. Set $v_0=a|z|^\delta$,
$V_0=\sigma\chi_{(-\eta,\eta)}(z)$, and $\mathscr V_0=v_0+V_0$, for $z\in(-h,h)$, with
$\sigma=(\sigma^1,\sigma^2)$, where $\chi_{(-\eta, \eta)}(z)$ is the
characteristic function of
the interval $(-h, h)$. Extend $v_0$ and $V_0$, and
consequently $\mathscr V_0$, periodically to the
whole space, and still use the same notations to denote the
extensions. Then, one can easily check that $v_0\in
C_\text{per}(\bar\Omega)\cap\mathcal H$ and $V_0\in
L^\infty(\Omega)\cap\mathcal H$.
By Corollary \ref{thm1.2},
there is a positive constant
$\varepsilon_0=\varepsilon_0(a,\delta,\eta,h)$, such that,
for any $\sigma=(\sigma^1,\sigma^2)$, with
$0<|\sigma|\leq\varepsilon_0$, system (\ref{5PESM})--(\ref{5ic})
has a unique weak solution, with initial data $\mathscr V_0$.
Note that
$$
\mathscr V_0=a|z|^\delta+\sigma\chi_{(-\eta,\eta)}(z), \quad z\in(-h,h).
$$
One can easily verify that $\mathscr V_0$ lies neither
in $X\cap\mathcal H$ nor
in $C_\text{per}(\bar\Omega)\cap\mathcal H$, and thus the results
established in \cite{BGMR03,TACHIM,PTZ09,KPRZ} cannot be applied to
prove the uniqueness of weak solutions with such kind of initial data.
\end{example}

\section{A tropical atmosphere model}
\label{SECTAMINT}
In the tropical zone of the earth, the wind in the lower troposphere is of equal magnitude, but with opposite sign to that in the upper troposphere, in other words, the primary effect is captured in the first baroclinic mode, i.e.\,on the first mode of the fluctuation of the solution about its vertical average. However, for the study of the tropical-extratropical interactions, where the transport of momentum between the barotropic (the spatial vertical average of the solution) and baroclinic (the fluctuation of the solution about the barotropic part) modes plays an important role, it is necessary to retain both the barotropic and baroclinic modes of the velocity.

Consider the following version of PEs (see, e.g., \cite{HAWI,LEWAN,MAJBOOK,PED,VALLIS,WP,ZENG}) for the atmosphere, in the layer $\mathbb R^2\times(0,H)$, for a positive constant $H$,
\begin{equation}\label{PE}
  \left\{
  \begin{array}{l}
    \partial_t\textbf{V}+(\textbf{V}\cdot\nabla_h)\textbf{V}+W\partial_z\textbf{V}
     -\mu\Delta\textbf{V}+\nabla_h\Phi=0,\\
    \partial_z\Phi=\frac{g\Theta}{\theta_0},\\
    \partial_t\Theta+\textbf{V}\cdot\nabla_h\Theta+W\partial_z\Theta +\frac{N^2\theta_0}{g}W=S_{\Theta},\\
    \nabla_h\cdot\textbf{V}+\partial_zW=0,
  \end{array}
  \right.
\end{equation}
where the unknowns $\textbf{V}=(V_1,V_2)^T$, $W$, $\Phi$ and $\Theta$ are the horizontal velocity field, vertical velocity, pressure and potential temperature, respectively, while the positive constant $\mu$ is the viscosity coefficient. The total potential temperature is given by
$$
\Theta^{\text{total}}(x,y,z,t)=\theta_0+\bar\theta(z)+\Theta(x,y,z,t),
$$
where $\theta_0$ is a positive reference constant temperature and $\bar\theta$ defines the vertical profile background stratification, satisfying $N^2=(g/\theta_0)\partial_z\bar\theta>0$, where $N$ is the Brunt-V\"ais\"al\"a buoyancy frequency. Here we use $\nabla_h=(\partial_x,\partial_y)$ to denote the horizontal gradient and $\textbf{V}^\perp=(-V_2, V_1)^T$.

Recall that the main effect of the tropical atmosphere is captured in the first baroclinic mode. By taking also the tropical-extratropical interactions into consideration, we can impose an ansatz of the form
\begin{equation*}
  \begin{pmatrix}
  \textbf{V} \\
  \Phi \\
\end{pmatrix}(x,y,z,t)=
\begin{pmatrix}
  u \\
  p \\
\end{pmatrix}(x,y,t)+
\begin{pmatrix}
  v \\
  p_1 \\
\end{pmatrix}(x,y,t)\sqrt2\cos(\pi z/H)
\end{equation*}
and
\begin{equation*}
  \begin{pmatrix}
  W \\
  \Theta \\
\end{pmatrix}(x,y,z,t)=
\begin{pmatrix}
  w \\
  \theta \\
\end{pmatrix}(x,y,t)\sqrt2\sin(\pi z/H),
\end{equation*}
which carry the barotropic and first baroclinic modes of the unknowns.

By performing the Galerkin projection of the PEs in the vertical direction onto the barotropic mode and the first baroclinic mode, one derives the following dimensionless interaction, between the barotropic mode and the first baroclinic mode, system  for the tropical atmosphere (see \cite{MAJBOOK}, and also \cite{MAJBIE,FRIMAJPAU,KHOMAJ1,STEMAJ}, for the details):
\begin{equation}\label{1.1}
\left\{
\begin{array}{l}
  \partial_tu+(u\cdot\nabla)u-\Delta u+\nabla p+\nabla\cdot(v\otimes v) =0,\\
  \nabla\cdot u=0,\\
  \partial_tv+(u\cdot\nabla)v-\Delta v+(v\cdot\nabla)u =\nabla\theta,\\
  \partial_t\theta+u\cdot\nabla\theta-\nabla\cdot v=S_\theta,
  \end{array}
  \right.
\end{equation}
where $u=(u^1, u^2)$ is the barotropic velocity, and $v=(v^1, v^2)$,
$p$ and $\theta$, respectively, are the first baroclinic modes of the
velocity, pressure and the temperature. The system is now defined on
$\mathbb R^2$, and the operators $\nabla$ and $\Delta$ are therefore
those for the variables $x$ and $y$.

An important ingredient of the tropical atmospheric circulation is the water vapour. Water vapour is the most abundant greenhouse gas in the atmosphere, and it is responsible for  amplifying the long-term warming or cooling cycles. Therefore, one should also consider the coupling with an equation modeling moisture in the atmosphere.

Following \cite{FRIMAJPAU}, we couple  system (\ref{1.1}) with the following large-scale moisture equation
\begin{equation}
  \partial_tq+u\cdot\nabla q+\bar Q\nabla\cdot v=-P,\label{1.5}
\end{equation}
where $\bar Q$ is the prescribed gross moisture stratification. The precipitation rate $P$ is parameterized, according to \cite{FRIMAJPAU,STEMAJ,NEEZEN,KHOMAJ2}, as
\begin{equation}\label{EQ}
P=\frac{1}{\varepsilon}(q-\alpha\theta-\hat q)^+,
\end{equation}
where $f^+=\max\{f, 0\}$ denotes the positive part of $f$, $\varepsilon$ is a convective adjustment time scale parameter, which is usually very small, and $\alpha$ and $\hat q$ are constants, with $\hat q>0$.

In order to close system (\ref{1.1})--(\ref{1.5}), one still needs to parameterize the source term $S_\theta$ in the temperature equation. Generally, the temperature source $S_\theta$ combines  three kinds of effects: the radiative cooling, the sensible heat flux and the precipitation $P$. For simplicity, and as in  \cite{FRIMAJPAU,MAJSOU}, we only consider in this paper the precipitation source term, i.e., we set
$$
S_\theta=P,
$$
with $P$ given by (\ref{EQ}).

As in \cite{FRIMAJPAU,MAJSOU}, by introducing the equivalent temperature $T_e$ and the equivalent moisture $q_e$ as
$$
T_e=q+\theta,\quad q_e=q-\alpha\theta-\hat q,
$$
system (\ref{1.1})--(\ref{1.5}) can be rewritten as
\begin{eqnarray}
  &&\partial_tu+(u\cdot\nabla)u-\Delta u+\nabla p+\nabla\cdot(v\otimes v)=0,\label{eq1}\\
  &&\nabla\cdot u=0,\label{eq2}\\
  &&\partial_tv+(u\cdot\nabla)v-\Delta v+(v\cdot\nabla)u=\frac{1}{1+\alpha}\nabla(T_e-q_e),\label{eq3}\\
  &&\partial_tT_e+u\cdot\nabla T_e-(1-\bar Q)\nabla\cdot v=0,\label{eq4}\\
  &&\partial_tq_e+u\cdot\nabla q_e+(\bar Q+\alpha)\nabla\cdot v=-\frac{1+\alpha}{\varepsilon} q_e^+, \label{eq5}
\end{eqnarray}
in $\mathbb R^2\times(0,\infty)$, where the constants $\alpha$ and $\bar Q$ are required to satisfy (see \cite{FRIMAJPAU})
\begin{equation}
  \label{req}
   0<\bar Q<1,\quad \alpha+\bar Q>0.
\end{equation}

Recalling that the adjustment relaxation time $\varepsilon$ in (\ref{eq5}) is small, both physically and mathematically, it is important to study the behavior of the solutions to system (\ref{eq1})--(\ref{eq5}), when $\varepsilon$ approaches to zero. Formally, by taking the relaxation limit, as $\varepsilon\rightarrow0^+$, system (\ref{eq1})--(\ref{eq5}) will converge to the following limiting system
\begin{eqnarray}
  &&\partial_tu+(u\cdot\nabla)u-\mu\Delta u+\nabla p+\nabla\cdot(v\otimes v)=0,\label{ineq1}\\
  &&\nabla\cdot u=0,\label{ineq2}\\
  &&\partial_tv+(u\cdot\nabla)v-\mu\Delta v+(v\cdot\nabla)u=\frac{1}{1+\alpha}\nabla(T_e-q_e),\label{ineq3}\\
  &&\partial_tT_e+u\cdot\nabla T_e-(1-\bar Q)\nabla\cdot v=0,\label{ineq4}\\
  &&\partial_tq_e+u\cdot\nabla q_e+(\bar Q+\alpha)\nabla\cdot v\leq 0,\label{ineq5}\\
  &&q_e\leq0,\label{ineq6}\\
  &&\partial_tq_e+u\cdot\nabla q_e+(\bar Q+\alpha)\nabla\cdot v= 0,\quad \mbox{a.e.~on }\{q_e<0\}. \label{ineq7}
\end{eqnarray}
Note that equation (\ref{eq5}) is now replaced by three inequalities (\ref{ineq5})--(\ref{ineq7}).

Inequality (\ref{ineq5}) comes from equation (\ref{eq5}), by noticing the negativity of the term $-\frac{1+\alpha}{\varepsilon}q_e^+$, while inequality
(\ref{ineq6}) is derived by multiplying both sides of equation (\ref{eq5}) by
$\varepsilon$, and taking the formal limit $\varepsilon\rightarrow0^+$. Inequality
(\ref{ineq7}) can be derived by the following heuristic argument: Let
$(u_\varepsilon, v_\varepsilon, T_{e\varepsilon}, q_{e\varepsilon})$ be a solution
to system (\ref{eq1})--(\ref{eq5}), and suppose that $(u_\varepsilon,
v_\varepsilon, T_{e\varepsilon}, q_{e\varepsilon})$ converges to $(u, v, T_e,
q_e)$, with $q_e\leq0$; for any compact subset $K$ of the set $\{(x,y,t)\in\mathbb R^2\times(0,\infty)~|~q_e(x,y,t)<0\}$, since
$q_{e\varepsilon}$ converges to $q_e$, one may have  $q_{e\varepsilon}<0$ on $K$,
for sufficiently small positive $\varepsilon$; therefore, by equation (\ref{eq5}),
it follows that $\partial_tq_{e\varepsilon}+u_\varepsilon\cdot\nabla
q_{e\varepsilon}+(\bar Q+\alpha)\nabla\cdot v_\varepsilon=0$, a.e.~on $K$, from
which, by taking $\varepsilon\rightarrow0^+$, one can see that (\ref{ineq7}) is
satisfied, a.e.~on $K$, and further a.e.~on $\{q_e<0\}$.

\section{Global well-posedness of a tropical atmosphere model}
\label{SECTAM}
In this section, we address the global well-posedness of strong solutions to the Cauchy problem of the tropical atmosphere model (\ref{eq1})--(\ref{eq5}) and its limiting system (\ref{ineq1})--(\ref{ineq7}). Strong convergence of the
relaxation limit from (\ref{eq1})--(\ref{eq5}) to the limiting system (\ref{ineq1})--(\ref{ineq7}) is also established.

\begin{definition}
  \label{def}
  Given a positive time $\mathcal T$ and the initial data $(u_0, v_0, T_{e,0}, q_{e,0})$. A function $(u, v, T_e, q_e)$ is called a strong solution to system (\ref{eq1})--(\ref{eq5}), on $\mathbb R^2\times(0,\mathcal T)$, with initial data $(u_0, v_0, T_{e,0}, q_{e,0})$, if it enjoys the following regularities
  \begin{eqnarray*}
  &(u, v)\in C([0,\mathcal T]; H^1(\mathbb R^1))\cap L^2(0,\mathcal T; H^2(\mathbb R^2)),\\
  &(\partial_tu,\partial_tv,\partial_tT_e, \partial_tq_e)\in L^2(0,\mathcal T; L^2(\mathbb R^2)), \\
  &(T_e,q_e)\in C([0,\mathcal T]; L^2(\mathbb R^2))\cap L^\infty(0,\mathcal T; H^1(\mathbb R^2)),
\end{eqnarray*}
and satisfies equations (\ref{eq1})--(\ref{eq5}), a.e.~on $\mathbb R^2\times(0,\mathcal T)$, and has the initial value
$$(u,v,T_e,q_e)|_{t=0}=(u_0, v_0, T_{e,0}, q_{e,0}).$$
\end{definition}

\begin{definition}
A function $(u, v, T_e, q_e)$ is called a global strong solution to system (\ref{eq1})--(\ref{eq5}), if it is a strong solution to system (\ref{eq1})--(\ref{eq5}), on $\mathbb R^2\times(0,\mathcal T)$, for any positive time $\mathcal T$.
\end{definition}

We have the following theorem on the global existence, uniqueness and well-posedness of strong solutions to the Cauchy problem of system (\ref{eq1})--(\ref{eq5}):

\begin{theorem}[cf. \cite{LITITITAMMOISTURE}]
  \label{glopositive}
Suppose that (\ref{req}) holds, and the initial data
\begin{equation}\label{asum}
(u_0, v_0, T_{e,0}, q_{e,0})\in H^1(\mathbb R^2),\quad\mbox{with}\quad \nabla\cdot u_0=0.
\end{equation}
Then, we have the following:

(i) There is a unique global strong solution $(u,v,T_e,q_e)$ to system (\ref{eq1})--(\ref{eq5}), with initial data $(u_0, v_0, T_{e,0}, q_{e,0})$, such that
\begin{align*}
\sup_{0\leq t\leq\mathcal T}&\|(u, v,T_e,q_e)(t)\|_{H^1}^2
+\int_0^\mathcal T\bigg(\frac{\|q_e^+\|_{H^1}^2}{\varepsilon}+
\|(u,v)\|_{H^2}^2+\|\nabla u\|_\infty\bigg)dt\\
&+\int_0^\mathcal T\|(\partial_tu,\partial_tv,\partial_t T_e)\|_2^2  dt
\leq C\left(\alpha,\bar Q, \mathcal T, \|(u_0, v_0, T_{e,0}, q_{e,0})\|_{H^1}\right),
\end{align*}
for any positive time $\mathcal T$, here and what follows, we use $C(\cdots)$ to denote a general positive constant depending only on the quantities in the parenthesis.

(ii) Suppose, in addition to (\ref{asum}), that $q_{e,0}\leq0$, a.e.~on $\mathbb R^2$, then
$$
\sup_{0\leq t\leq \mathcal T}\frac{\|q_e^+(t)\|_2^2}{\varepsilon} +\int_0^\mathcal T\|\partial_tq_e\|_2^2dt
\leq C\left(\alpha,\bar Q, \mathcal T, \|(u_0, v_0, T_{e,0}, q_{e,0})\|_{H^1}\right),
$$
for any positive time $\mathcal T$.

(iii) Suppose, in addition to (\ref{asum}), that $(\nabla T_{e,0},\nabla q_{e,0})\in L^m(\mathbb R^2)$, for some $m\in(2,\infty)$, then the following estimate holds
\begin{align*}
  \sup_{0\leq t\leq\mathcal T}\|(\nabla T_{e},\nabla q_{e})(t)\|_m^2
  \leq C\left(\alpha,\bar Q, \mathcal T, \|(u_0, v_0, T_{e,0}, q_{e,0})\|_{H^1},\|(\nabla T_{e,0},\nabla q_{e,0})\|_m\right),
\end{align*}
for any positive time $\mathcal T$, and the unique strong solution $(u,v, T_e, q_e)$ depends continuously on the initial data, on any finite interval of time.
\end{theorem}

Global existence of strong solutions are based on those a priori
estimates stated in (i) of Theorem \ref{glopositive}. We successfully
do the $L^\infty(L^2)$ estimate on $(u, v, T_e, q_e)$, $L^\infty(L^4)$
estimate on $(u, v, T_e, q_e)$, $L^\infty(H^1)$ estimate on $u$, and
finally the $L^\infty(H^1)$ estimate on $(u, v, T_e, q_e)$, and
each of them bases on the previous ones.

Next, we focus on the uniqueness. Let $(u, v, T_e, q_e)$ and $(\tilde u, \tilde v, \tilde T_e, \tilde q_e)$ be two strong solutions to system (\ref{eq1})--(\ref{eq5}), with the same initial data $(u_0, v_0, T_{e,0}, q_{e, 0})$, on the time interval $(0,\mathcal T)$. Define the new functions
$$
(\delta u, \delta v, \delta T_e,\delta q_e)=(u, v, T_e, q_e)-(\tilde u, \tilde v, \tilde T_e, \tilde q_e).
$$
By studying the system for $(\delta u, \delta v, \delta T_e,\delta q_e)$, we have
\begin{eqnarray}
&&\frac{d}{dt}\|(\delta u,\delta v,\delta T_e,\delta q_e)\|_2^2+\frac14\|(\delta u,\delta v)\|_{H^1}^2\nonumber\\
&\leq& C\left(1+\|(\tilde u, \tilde v)\|_4^4+\|(\nabla\tilde u, \nabla\tilde v,\nabla v)\|_2^2\right)\|(\delta u,\delta v,\delta T_e, \delta q_e)\|_2^2\nonumber\\
&&+C\|(\nabla\tilde T_e,\nabla\tilde q_e)\|_2\|(\delta u,\delta v)\|_\infty\|(\delta T_e, \delta q_e)\|_2.\label{LINEW6.2}
\end{eqnarray}
Recalling the Brezis--Gallouet--Wainger inequality (see \cite{Brezis_Gallouet_1980,Brezis_Wainger_1980})
$$
\|f\|_{L^\infty(\mathbb R^2)}\leq C\|f\|_{H^1(\mathbb R^2)}\log^{\frac12}\left(\frac{\|f\|_{H^2(\mathbb R^2)}}{\|f\|_{H^1(\mathbb R^2)}}+e\right),
$$
and denoting $U=(u,v), \tilde U=(\tilde u, \tilde v)$ and $\delta U=(\delta u, \delta v)$, we have
\begin{eqnarray}
\|\delta U\|_\infty&\leq& C\left[\|\delta U\|_{H^1}^2\log^+\left(\frac{S(t)}{\|\delta U\|_{H^1}}\right)\right]^{\frac12}, \label{LINEW6.2-1}
\end{eqnarray}
where
$$
S(t)=\|U\|_{H^2}+\|\tilde U\|_{H^2}+e(\|U\|_{H^1}+\|\tilde U\|_{H^1}).
$$
Denoting
\begin{eqnarray*}
  &&f=\|(\delta u,\delta v,\delta T_e,\delta q_e)\|_2^2,\quad G=\frac14\|(\delta u,\delta v)\|_{H^1}^2,\\
  &&m_1=C\left(1+\|(\tilde u, \tilde v)\|_4^4+\|(\nabla\tilde u, \nabla\tilde v,\nabla v)\|_2^2\right),\quad m_2=C\|(\nabla\tilde T_e,\nabla\tilde q_e)\|_2,
\end{eqnarray*}
then it follows from (\ref{LINEW6.2}) and (\ref{LINEW6.2-1}) that
$$
f'+G\leq m_1f+m_2\left[fG\log^+\left(\frac{S/4}{G}\right)\right]^{\frac12}.
$$

Then, the uniqueness follows from the next lemma.

\begin{lemma}[cf. \cite{LITITITAMMOISTURE}]\label{lemuniq}
Given a positive time $\mathcal T$, and let $m_1, m_2$ and $S$ be nonnegative functions on $(0,\mathcal T)$, such that
$$
m_1, S\in L^1((0,\mathcal T)),\quad m_2\in L^2((0,\mathcal T)), \mbox{ and } S>0, \mbox{ a.e.~on } (0,\mathcal T).
$$
Suppose that $f$ and $G$ are two nonnegative functions on $(0,\mathcal T)$, with $f$ being absolutely continuous on $[0,\mathcal T)$, and satisfy
\begin{equation*}
\left\{
\begin{array}{l}
f'(t)+G(t)\leq m_1(t)f(t)+m_2(t)\left[f(t)G(t)\log^+\left(\frac{S(t)}{G(t)}\right) \right]^{\frac12},\quad\mbox{ a.e.~on }(0,\mathcal T),\\
f(0)=0,
\end{array}
\right.
\end{equation*}
where $\log^+z=\max\{0,\log z\}$, for $z\in(0,\infty)$, and when $G(t)=0$, at some time $t\in[0,\mathcal T)$, we adopt the following natural convention
$$
G(t)\log^+\left(\frac{S(t)}{G(t)}\right)=\lim_{z\rightarrow 0^+}z\log^+\left(\frac{S(t)}{z}\right)=0.
$$
Then, we have $f\equiv0$ on $[0,\mathcal T)$.
\end{lemma}

We also have the global existence and uniqueness of strong solutions
to the limiting system (\ref{ineq1})--(\ref{ineq7}), where
strong solutions
to system (\ref{ineq1})--(\ref{ineq7}) are defined in the similar way
as those to system (\ref{eq1})--(\ref{eq5}). In fact, we have the following theorem:

\begin{theorem}[cf. \cite{LITITITAMMOISTURE}]\label{glozero}
Suppose that (\ref{req}) holds, and the initial data
\begin{equation}\label{asum1}
(u_0, v_0, T_{e,0}, q_{e,0})\in H^1(\mathbb R^2),\quad \nabla\cdot u_0=0,\quad q_{e,0}\leq0,\mbox{ a.e.~on }\mathbb R^2.
\end{equation}
Then, there is a unique global strong solution $(u,v,T_e,q_e)$ to system (\ref{ineq1})--(\ref{ineq7}), with initial data $(u_0, v_0, T_{e,0}, q_{e,0})$, such that
\begin{align*}
  \sup_{0\leq t\leq\mathcal T} \|(u, v,T_{e},q_{e})(t)\|_{H^1}^2
  & +\int_0^\mathcal T\left(
  \|(u,v)\|_{H^2}^2+\|\nabla u\|_\infty+\|(\partial_tu,\partial_tv,\partial_tT_e,
  \partial_tq_{e})\|_2^2\right) dt\\
  \leq&C\left(\alpha,\bar Q, \mathcal T, \|(u_0, v_0, T_{e,0}, q_{e,0})\|_{H^1}\right),
\end{align*}
for any positive time $\mathcal T$.

If we assume, in addition, that $(\nabla T_{e,0}, \nabla q_{e,0})\in L^m(\mathbb R^2)$, for some $m\in(2,\infty)$, then we have further that
\begin{align*}
  \sup_{0\leq t\leq\mathcal T}\|(\nabla T_{e},\nabla q_{e})(t)\|_m^2
  \leq C\left(\alpha,\bar Q, \mathcal T, \|(u_0, v_0, T_{e,0}, q_{e,0})\|_{H^1},\|(\nabla T_{e,0},\nabla q_{e,0})\|_m\right),
\end{align*}
for any positive time $\mathcal T$, and the unique strong solution $(u, v, T_e, q_e)$ depends continuously on the initial data.
\end{theorem}

The existence part of Theorem \ref{glozero} is proven by taking the limit to strong solutions of (\ref{eq1})--(\ref{eq5}). Denote by $(u_\varepsilon, v_\varepsilon, T_{e\varepsilon}, q_{e\varepsilon})$ the unique strong solution to system (\ref{eq1})--(\ref{eq5}), with initial data $(u_0, v_0, T_{e0}, q_{e0})$. Thanks to the a priori estimates stated in Theorem \ref{glopositive}, the strong solutions $(u_\varepsilon, v_\varepsilon, T_{e\varepsilon}, q_{e\varepsilon})$
converge to some $(u, v, T_{e}, q_{e})$. Thanks to the a priori estimates stated in Theorem \ref{glopositive}, by using Aubin-Lions compactness lemma (see, e.g., Simon \cite{Simon}), one can show that $(u, v, T_e, q_e)$ satisfies (\ref{ineq1})--(\ref{ineq6}). It remains to verify (\ref{ineq7}). To this end, let us define the set
$$
\mathcal O^-=\{(x,t)|q_e(x,t)<0, x\in\mathbb R^2, t\in(0,\infty)\},
$$
and for any positive integers $j,k,l$, we define
$$
\mathcal O^-_{jkl}=\left\{(x,t)\bigg|q_e(x,t)<-\frac1j, x\in B_k, t\in(0,l)\right\},
$$
where $B_k \subset \mathbb{R}^2$ is a disc of radius $k$, and  $j,k,l \in \mathbb{N}$.
Noticing that
$$
\mathcal O^-=\cup_{j}^\infty\cup_{k=1}^\infty\cup_{l=1}^\infty \mathcal O^-_{jkl},
$$
to prove that (\ref{ineq7}) holds a.e.~on $\mathcal O^-$, it suffices
to show that it holds a.e.~on $\mathcal O^-_{jkl}$, for any positive
integers $j,k,l$.
Now, let us fix the positive integers $j,k,l$. Thanks to the a priori estimate in Theorem \ref{glopositive}, by the Aubin-Lions compactness lemma (see, e.g., Simon \cite{Simon}), and using the Cantor diagonal argument, one can show that
$q_{e\varepsilon}\rightarrow q_e$ in $C([0,\mathcal T]; L^2(B_R))$,
for any positive time $\mathcal T$ and positive radius $R$, and thus $q_{e\varepsilon}\rightarrow q_e$ in
$L^2(\Omega_{jkl})$. Therefore, there is a subsequence, still denoted
by $q_{e\varepsilon}$, such that $q_{e\varepsilon}\rightarrow q_e$,
a.e.~on $\mathcal O^-_{jkl}$. By the Egoroff theorem, for any positive
number $\eta>0$, there is a subset $E_\eta$ of $\mathcal O^-_{jkl}$,
with $|E_\eta|\leq\eta$, such that
$$
q_{e\varepsilon}\rightarrow q_e,\quad\mbox{uniformly on }\mathcal O^-_{jkl}\setminus E_\eta.
$$
Recalling the definition of $\mathcal O^-_{jkl}$, this implies that for sufficiently small positive $\varepsilon$, it holds that
$$
q_{e\varepsilon}\leq q_e+\frac{1}{2j}\leq-\frac{1}{2j}<0, \quad \mbox{on }\mathcal O^-_{jkl}\setminus E_\eta.
$$
As a result, by equation (\ref{eq5}) for $q_{e\varepsilon}$, we have, for any sufficiently small positive $\varepsilon$, that
$$
\mathcal G_\varepsilon:=\partial_tq_{e\varepsilon}+u_\varepsilon\cdot\nabla q_{e\varepsilon}+(\bar Q+\alpha)\nabla\cdot v_\varepsilon=0, \quad\mbox{ a.e.~on }\mathcal O^-_{jkl}\setminus E_\eta.
$$
Noticing that
$$
\mathcal G_\varepsilon\rightharpoonup  \partial_tq_{e }+u \cdot\nabla q_{e }+(\bar Q+\alpha)\nabla\cdot v=:\mathcal G ,\quad\mbox{ in }L^2(0,\mathcal T; L^2(\mathbb R^2)),
$$
for any positive finite time $\mathcal T$, which in particular implies $\mathcal G_\varepsilon\rightharpoonup\mathcal G$, in $L^2(\mathcal O_{jkl}\setminus E_\eta)$. Since $\mathcal G_\varepsilon=0$, a.e.~on $\mathcal O_{jkl}\setminus E_\eta$, we have $\mathcal G=0$, a.e.~on $\mathcal O_{jkl}\setminus E_\eta$, that is
$$
\partial_tq_{e }+u \cdot\nabla q_{e }+(\bar Q+\alpha)\nabla\cdot v=0, \quad\mbox{ a.e.~on }\Omega_{jkl}\setminus E_\eta.
$$
This implies that the above equation holds, a.e.~on $\mathcal O^-_{jkl}$, and further on $\mathcal O^-$, in other words, (\ref{ineq7}) holds.

For the uniqueness part, some additional attention is needed to be paid to $q_e$, because it satisfies equation (\ref{ineq7}) only on part of the domain, i.e.\,on $\{q_e<0\}$. Let $(u, v, T_e, q_e)$ and $(\tilde u, \tilde v, \tilde T_e, \tilde q_e)$ be two strong solutions to system (\ref{ineq1})--(\ref{ineq7}), with the same initial data $(u_0, v_0, T_{e,0}, q_{e, 0})$. Define the new functions
$$
(\delta u, \delta v, \delta T_e,\delta q_e)=(u, v, T_e, q_e)-(\tilde u, \tilde v, \tilde T_e, \tilde q_e).
$$
Performing the energy estimate to the system for $(\delta u,\delta v,\delta T_e)$, one obtains
\begin{align}
  &\frac{d}{dt}\|(\delta u,\delta v,\delta T_e)\|_2^2+\|\nabla\delta u\|_2^2+\|\nabla\delta v\|_2^2\nonumber\\
  \leq&C\int_{\mathbb R^2}[(|\nabla\tilde u|+|\nabla\tilde v|+|\nabla v|+|v|^2+|\tilde v|^2)(|\delta u|^2+|\delta v|^2)\nonumber\\
  &+|\delta T_e|^2+|\delta q_e|^2+|\nabla\tilde T_e||\delta u||\delta T_e|]dxdy.\label{uni1'}
\end{align}
To derive the equation for $\delta q_e$, let us recall the well-known fact that the derivatives of a function $f\in W^{1,1}_{\text{loc}}(\Omega)$ vanish, a.e.~on any level set $\{(x,y,t)\in\Omega|f(x,y,t)=c\}$, see, e.g., \cite{EVANS} or page 297 of \cite{GIOVANNI}. By the aid of this fact, and using equation (\ref{ineq7}), one obtains the equation for $\delta q_e$ as
\begin{eqnarray}
&\partial_t\delta q_e+u\cdot\nabla\delta q_e+\delta u\cdot\nabla\tilde q_e
  = -(\bar Q+\alpha)[\nabla\cdot\delta v\chi_{\Omega_1}
  +\nabla\cdot v\chi_{\Omega_2}-\nabla\cdot\tilde v\chi_{\Omega_3}]\nonumber\\
  &=-(\bar Q+\alpha)[\nabla\cdot\delta v-\nabla\cdot\delta v\chi_{\Omega_4}
  +\nabla\cdot \tilde v\chi_{\Omega_2}-\nabla\cdot v\chi_{\Omega_3}],\label{deq5}
\end{eqnarray}
a.e.~on $\Omega=\mathbb R^2\times(0,\infty)$, where
\begin{eqnarray*}
  &\Omega_1=\{q_e<0\}\cap\{\tilde q_e<0\}, \quad \Omega_2=\{q_e<0\}\cap\{\tilde q_e=0\},\\
  &\Omega_3=\{q_e=0\}\cap\{\tilde q_e<0\},\quad \Omega_4=\{q_e=0\}\cap\{\tilde q_e=0\}.
\end{eqnarray*}
Performing $L^2$ energy estimate to (\ref{deq5}), one obtains
\begin{align}
  \frac12\frac{d}{dt}\|\delta q_e\|_2^2
  \leq&\frac14\int_{\mathbb R^2}|\nabla\delta v|^2dxdy+C\int_{\mathbb R^2}(|\delta q_e|^2+|\nabla\tilde q_e||\delta u||\delta q_e|)dxdy\nonumber\\
  &-(\bar Q+\alpha)\int_{\mathbb R^2}(\nabla\cdot\tilde v\chi_{\Omega_2}-\nabla\cdot v\chi_{\Omega_3})(q_e-\tilde q_e)dxdy. \label{uni2}
\end{align}
Note that, in the last integral on the right-hand side of the above inequality, it is linear in $\delta q_e$, while other terms in (\ref{uni2}) are quadratic in $\delta q_e$, and as a result, the above inequality, combined with (\ref{uni1'}), does not necessary imply, in general, the uniqueness. Fortunately, thanks to inequality (\ref{ineq5}), one can actually show that the last integral on the right-hand side of the above inequality is of positive sign, and therefore it can be ignored there, in other words, we have
\begin{align}
  \frac12\frac{d}{dt}\|\delta q_e\|_2^2
  \leq&\frac14\int_{\mathbb R^2}|\nabla\delta v|^2dxdy+C\int_{\mathbb R^2}(|\delta q_e|^2+|\nabla\tilde q_e||\delta u||\delta q_e|)dxdy. \label{uni2'}
\end{align}
Combining (\ref{uni1'}) and (\ref{uni2'}), one obtains again (\ref{LINEW6.2}), and thus proves the uniqueness.

While the strong convergence of the relaxation limit, as $\varepsilon\rightarrow0$, of system (\ref{eq1})--(\ref{eq5}) to the limiting system (\ref{ineq1})--(\ref{ineq7}) is stated in the next theorem:

\begin{theorem}[cf. \cite{LITITITAMMOISTURE}]
  \label{cong}
 Suppose that (\ref{req}) holds and the initial data
\begin{eqnarray*}
  &&(u_0, v_0, T_{e,0}, q_{e,0})\in H^1(\mathbb R^2),\quad \nabla\cdot u_0=0,\\
  &&(\nabla T_{e,0},\nabla q_{e,0})\in L^m(\mathbb R^2),\quad q_{e,0}\leq0,\mbox{ a.e.~on }\mathbb R^2,
\end{eqnarray*}
for some $m\in(2,\infty)$. Denote by $(u_\varepsilon, v_\varepsilon, T_{e\varepsilon}, q_{e\varepsilon})$ and $(u, v, T_e, q_e)$ the unique global strong solutions to systems (\ref{eq1})--(\ref{eq5}) and (\ref{ineq1})--(\ref{ineq7}), respectively, with the same initial data $(u_0, v_0, T_{e,0}, q_{e,0})$.

Then, we have the estimate
\begin{align*}
  &\sup_{0\leq t\leq\mathcal T}\|(u_\varepsilon-u, v_\varepsilon-v, T_{e\varepsilon}-T_e,q_{e\varepsilon}-q_e)(t)\|_2^2\\
  &+\int_0^\mathcal T\left(\|(\nabla(u_\varepsilon-u),\nabla(v_\varepsilon-v))\|_2^2 +\frac{\|q_{e\varepsilon}^+\|_2^2}{\varepsilon}\right)dt\leq C\varepsilon,
\end{align*}
for any finite positive time $\mathcal T$, where $C$ is a positive constant depending only on $\alpha, \bar Q, m, \mathcal T$, and the initial norm $\|(u_0, v_0,q_{e,0}, T_{e,0})\|_{H^1}+\|(\nabla T_{e,0},\nabla q_{e,0})\|_m$.

Therefore, in particular, we have the convergences
\begin{eqnarray*}
  &(u_\varepsilon, v_\varepsilon)\rightarrow(u, v)\quad\mbox{in }L^\infty(0,\mathcal T; L^2(\mathbb R^2))\cap L^2(0,\mathcal T; H^1(\mathbb R^2)),\\
  &(T_{e\varepsilon}, q_{e\varepsilon})\rightarrow(T_e,q_e)\quad\mbox{in }L^\infty(0,\mathcal T; L^2(\mathbb R^2)),\quad q_{e\varepsilon}^+\rightarrow0\quad\mbox{in }L^2(0,\mathcal T; L^2(\mathbb R^2)),
\end{eqnarray*}
for any positive time $\mathcal T$, and the convergence rate is of order $O(\sqrt\varepsilon)$.
\end{theorem}

Define the difference function $(\delta u_\varepsilon, \delta v_\varepsilon, \delta T_{e\varepsilon}, \delta q_{e\varepsilon})$ as
$$
(\delta u_\varepsilon, \delta v_\varepsilon, \delta T_{e\varepsilon}, \delta q_{e\varepsilon})=(u_\varepsilon,v_\varepsilon, T_{e\varepsilon}, q_{e\varepsilon})-(u, v, T_e, q_e).
$$
Performing the energy estimates to the system for $(\delta u_\varepsilon, \delta v_\varepsilon, \delta T_{e\varepsilon})$, one obtains
\begin{eqnarray}
  && \frac{d}{dt}\|(\delta u_\varepsilon,\delta v_\varepsilon,\delta T_{e\varepsilon})\|_2^2+\|(\nabla\delta u_\varepsilon,\nabla\delta v_\varepsilon)\|_2^2\nonumber\\
  &\leq&C\int_{\mathbb R^2}[(|\nabla u|+|\nabla v|+|v|^2)(|\delta u_\varepsilon|^2+|\delta v_\varepsilon|^2)\nonumber\\
  && +|\delta T_{e\varepsilon}|^2+|\delta q_{e\varepsilon}|^2+|\nabla T_e||\delta u_\varepsilon||\delta T_{e\varepsilon}|]dxdy.\label{con1}
\end{eqnarray}
As before, by the aid of the fact that the derivatives of a function are zero a.e.\,on any level set, and using equation (\ref{eq5}) for $q_{e\varepsilon}$ and equation (\ref{ineq7}) for $q_e$, one can derive the equation for $\delta q_{e\varepsilon}$ as
\begin{align*}
  \partial_t\delta q_{e\varepsilon}&+\delta u_\varepsilon\cdot\nabla\delta q_{e\varepsilon}+\delta u_\varepsilon\cdot\nabla q_e+u\cdot\nabla\delta q_{e\varepsilon}\nonumber\\
  &+(\bar Q+\alpha)\nabla\cdot\delta v_\varepsilon=-\frac{1+\alpha}{\varepsilon}q_{e\varepsilon}^+-(\bar Q+\alpha)\nabla\cdot v\chi_{\mathcal O}(x,y,t),
\end{align*}
a.e.~on $\mathbb R^2\times(0,\infty)$, where $\mathcal O=\{(x,t)\in\mathbb R^2\times(0,\infty)|q_e(x,t)=0\}$. Performing the $L^2$ energy estimates to the above system, and using same similar tricks as before, one can obtain
\begin{align}
  \frac{d}{dt}\|\delta q_{e\varepsilon}\|_2^2+\frac{1+\alpha}{\varepsilon}\|q_{e\varepsilon}^+\|_2^2
  \leq&\frac12\|\nabla\delta v_\varepsilon\|_2^2+2(\bar Q+\alpha)^2\|\delta q_{e\varepsilon}\|_2^2+\frac{(\bar Q+\alpha)^2}{1+\alpha}\varepsilon\|\nabla v\|_2^2\nonumber\\
  &+2\int_{\mathbb R^2}|\nabla q_e||\delta u_\varepsilon||\delta q_{e\varepsilon}|dxdy.\label{con4}
\end{align}
Starting from (\ref{con1}) and (\ref{con4}), one can then obtain the desired estimate in Theorem \ref{cong}.

\begin{remark}
(i) In the absence of the barotropic mode, global existence and uniqueness of strong solutions to the inviscid limiting system was proved in \cite{MAJSOU}, and the relaxation limit, as $\varepsilon\rightarrow0^+$, was also studied there, but the convergence rate was not achieved. Note that in the absence of the barotropic mode, the limiting system is linear, while in the presence of the barotropic mode, the limiting system is nonlinear.

(ii) Existence and uniqueness of solutions to the limiting system (\ref{ineq1})--(\ref{ineq7}), without viscosity, was proposed as an open problem in \cite{FRIMAJPAU}, and also in \cite{MAJSOU,KHOMAJSTE,MAJAB}. Notably, Theorem \ref{glozero} settles this open problem for  the \textsc{viscous version} of (\ref{ineq1})--(\ref{ineq7}). Note that we only add viscosity to the velocity equations, and we do not use any diffusivity in the temperature and moisture equations.

(iii) For the case without the coupling with the moisture equations, i.e.\,for system (\ref{1.1}), global well-posedness of strong solutions was established in Li--Titi \cite{LITITITCM} by a different method
from that presented here. Global well-posedness of strong solutions to a coupled system of the primitive equations with moisture (therefore, it is a different system from those considered in this paper) was recently addressed in Coti Zelati et al \cite{ZHKTZ}, where the system under consideration has full dissipation in all dynamical equations, and in particular has diffusivity in the temperature and moisture equations.
\end{remark}

\section*{Acknowledgments}
{J.L. is thankful to the kind hospitality of Texas A\&M University where part of this work was completed. This work was supported in part by the ONR grant N00014-15-1-2333 and the NSF grants DMS-1109640 and DMS-1109645.}
\par

\end{document}